\documentclass[12pt]{article}
\usepackage{amsmath,amssymb,amscd,verbatim,amsthm}
\newtheorem{theorem}{Theorem}[section]

\theoremstyle{remark}

\theoremstyle{definition}

\def\Dbar{\leavevmode\lower.6ex\hbox to 0pt{\hskip-.23ex\accent"16\hss}D}

\newcommand\dom{\operatorname{dom}}

\newcommand\cf{\operatorname{cf}}
\newcommand\cof{\operatorname{cof}}
\newcommand\otp{\operatorname{otp}}

\def\Dbar{\leavevmode\lower.6ex\hbox to 0pt{\hskip-.23ex \accent"16\hss}D}
\begin{document}
\title{In memoriam: James Earl Baumgartner (1943--2011)}
\author{
J.A.~Larson
{Department of Mathematics}\\
{University of Florida, Gainesville}\\
{Gainesville, FL 32611--8105, USA}\\
}
\maketitle
\begin{abstract}
James Earl Baumgartner (March 23, 1943  -- December 28, 2011)
came of age mathematically during the emergence of forcing as a fundamental
technique of set theory, and his seminal research changed the way set
theory is done. He made fundamental contributions to the development
of forcing, to our understanding of uncountable orders, to the
partition calculus, and to large cardinals and their ideals.
He promulgated the use of logic such as absoluteness and elementary
submodels to solve problems in set theory, he applied his knowledge of set theory to a
variety of areas in collaboration with other mathematicians, and he
encouraged a community of mathematicians with engaging survey talks,
enthusiastic discussions of open problems, and friendly mathematical conversations.
\end{abstract}

\section{Overview of Baumgartner's Life}
James E. Baumgartner was born on March 23, 1943 in Wichita, Kansas.
His high school days included tennis, football, and leading roles in school plays.
In 1960 he entered the California Institute of Technology, but stayed only
two years, moving to the University of California, Berkeley in 1962,
in part because it was co-educational. There he met and married his
wife Yolanda.  He continued his interest in
drama and mathematics as an undergraduate, earned his A.B. 
in mathematics in 1964, and continued study as a graduate student. Baumgartner  \cite[page 2]{JEB1970thesis} 
dated his interest in set theory to the four week long 1967 UCLA Summer Institute on Axiomatic Set Theory.\footnote{John W. 
Addison, Jr.,who directed Baumgartner's earliest research, made it possible 
for Baumgartner to attend the 1967 UCLA Summer School.} 
The mathematics for his dissertation was completed in spring 1969, and  
Baumgartner became a John Wesley Young Instructor at Dartmouth College in 
fall 1969.  His adviser, Robert Vaught, required Baumgartner to include 
additional details in his dissertation, so he did not earn his doctorate 
until 1970.

In Fall 1970, he was a
Participating Scholar in the New York Academy of Sciences
Scholar-in-Residence Program under the sponsorship of Paul Erd\H{o}s.
In Fall 1971, his status at Dartmouth shifted from John Wesley Young
Instructor to Assistant Professor but he spent the academic year 1971-1972 
at the California Institute of Technology as a Visiting Assistant
Professor and regularly attended the UCLA Logic Colloquia on Fridays.
Winter and spring quarters of 1975 were spent at the University of
California, Berkeley as a Research Associate.  He was tenured and
promoted to Associate Professor in 1976, promoted to Professor in
1980, became the first John G. Kemeny Professor of Mathematics in
1983. Baumgartner, along with Donald~A. Martin, and Saharon Shelah,
organized the 1983 American Mathematical Society Summer Research 
Conference on \emph{Axiomatic Set Theory} in Boulder, Colorado,
and they edited the proceedings  \cite{BaxiomaticSetTheory1984}.
This event was the first large meeting devoted entirely to set theory 
since the 1967 Summer School held at UCLA.  
Baumgartner spent a stint as Chair of the department from 1995-1998.
He was honored with the Baumgartner Fest in 2003 at which a number
of his students spoke.  Slowed by the multiple sclerosis diagnosed 
in 1982, he retired with emeritus status in 2004. He died in 2011 
under the care of his loving wife, Yolanda. 

Baumgartner had ten doctoral students at Dartmouth listed below 
with academic affiliations for those who have one: 
Robert Beaudoin (1985), 
Stefan Bilaniuk (1989), Trent University, 
Denis Devlin (1980), 
Claudia Henrion (1985), 
Albin Jones (1999), 
Jean Larson (1972), University of Florida, 
Thomas Leathrum (1993), Jacksonville State University, Alabama,
Tadatoshi Miyamoto (1988), Nanzan University of Nagoya,
Alan Taylor (1975), Union College,
Stanley Wagon (1975), Macalester College.

Charles K. Landraitis, who is affiliated with Boston College, 
was a set theory student at Dartmouth College graduating in 1975,
and was often included in group activities for the Baumgartner
group of students, even though his adviser was Victor Harnik of Haifa 
University, who visited Dartmouth College.
Baumgartner also advised Peter Dordal  who graduated in 1982 from Harvard
University and is now in computer science at Loyola University in Chicago.

Baumgartner enjoyed working with a number of mathematicians on post-doctoral 
positions or visiting positions at Dartmouth, most as John Wesley
Young Instructors. These included
Uri Abraham, 
Alessandro Andretta, 
J\"org Brendle, 
Elizabeth Theta Brown, 
James Cummings, 
Frantisek Franek, 
Jean-Pierre Levinski, 
George McNulty, 
Lee Stanley, 
Claude Sureson, 
Stevo Todorcevic, 
Robert Van Wesep, and 
Jindrich Zapletal.  

\section{A personal note}
I interacted most with Baumgartner as a graduate student.
He arrived at Dartmouth College my second year in graduate school, 
and I think I met him in my oral qualifying exam where he clarified
a question I was being asked enabling me to answer it successfully.

I took his course in set theory starting that fall and decided to ask 
to work with him. Alas, I was the   second to ask, and since the first to ask 
quickly switched to someone else, at Baumgartner's suggestion, we 
informally started reading together.

Baumgartner had learned a lot from fellow students in graduate school,
so encouraged me to bring an undergraduate
into our  fall 1970 conversations, enabling me to have a peer with whom to talk. 
My fellow student took a term of independent study with Baumgartner
that fall, but after one semester, continued without credit since
he had spent more time than he felt he could afford on it.

In spring 1971 I came up with a short proof of  
$\omega^\omega\rightarrow (\omega^\omega,n)^2$ which became the 
cornerstone of my thesis 
even though 
C.C.~Chang had already proved it 
for $n=3$ and Eric Milner showed how to generalize his proof for all $n$.  
I had difficulty explaining the proof to Baumgartner, so week after week, 
he would tell me he did not yet understand, that he was sure I would be 
able to explain it to him, and he would cheerfully ask me to come back 
next week to try again.

Once he understood my proof, Baumgartner arranged support for me to 
attend the 1971 Summer School in Cambridge, organized by Adrian Mathias, where 
I met for the first time a very large number of mathematicians that I have
continued to see.  Baumgartner sent me with a paper of his to hand deliver 
to Hajnal, guaranteeing a meeting with Erd\H{o}s, Hajnal and Milner. 

My final year in graduate school was spent as a visiting graduate student 
at UCLA where Baumgartner shared his UCLA office with me 1971-72, 
while he spent most of his time at Cal Tech.

It was always wonderful to visit the Baumgartners in 
Hanover and at many conferences over the years.  
In October 2003, Arthur Apter and Marcia Groszek organized a 
Baumgartner Fest in honor of his 60th birthday. 
It was a wonderful conference with many people speaking on
mathematics of interest to Baumgartner.
At a party at his house during this conference, Baumgartner passed around the 
framed conference photo from the 1967 UCLA Summer Institute on Axiomatic Set 
Theory that was the beginning of his interest in set theory. 
It was a time to reflect on all the meetings we had enjoyed in 
between the 1967 UCLA Summer Institute and the Baumgartner Fest.

\section{Baumgartner's mathematical work}
We now turn to the mathematical context in which Baumgartner worked
and a discussion of a selected works mainly by date of publication.
 Jech's book \cite{Jech2003} has been generally followed for definitions and
 notation. Kanamori's book \cite{Kanamori2003} has been an invaluable
resource for both mathematics and history.

\subsection{Mathematical context and graduate school days}
Baumgartner \cite[462]{JEB1986onKunen} described the
mathematical scene in the years just prior to his
time as a set theory graduate student:
\begin{quote}
Once upon a time, not so very long ago, logicians hardly ever wrote anything down. Wonderful results were being obtained almost weekly, and no one wanted to miss out on the next theorem by spending the time to write up the last one. Fortunately there was a Center where these results were collected and organized, but even for the graduate students at the Center life was hard. They had no textbooks for elementary courses, and for advanced courses they were forced to rely on handwritten proof outlines, which were usually illegible and incomplete; handwritten seminar notes, which were usually wrong; and Ph.D. dissertations, which were usually out of date. Nevertheless, they prospered.
Now the Center I have in mind was Berkeley and the time was the early
and middle 1960's, \dots 
\end{quote}
In the early and middle 1960's aspects of set theory were 
developing in concert: forcing,\footnote{Forcing is a technique for adjoining 
a generic object to a given model of set theory so that the properties 
of the generic object and hence the extension of the original model 
generated from the generic object are determined by the construction 
and the original model.}  large cardinals, combinatorial set theory 
and interactions with model theory. 

Forcing was introduced by Paul Cohen in 1963,  and it 
was quickly applied by Easton to the question of the
size of powers of regular cardinals in his 1964 thesis \cite{Easton1964}.
Robert Solovay \cite{Solovay1964Julabs}, \cite{Solovay1970} proved
the consistency of ZF with every set of reals being Lebesgue
measurable by 1964, but only published the result in 1970.  
Early lecture series around the world included Prikry's January 1964
lecture in Mostowski's seminar in Warsaw, Levy's course on forcing in
1964, \cite[161]{MooreG1988} and lectures by Jensen at the University
of Bonn in 1965-66 \cite{Jensen1967book}.  

Large cardinal concepts date back to Hausdorff \cite{Hausdorff1908}
(weakly inaccessible),
Mahlo (Mahlo cardinals\footnote{$\kappa$ is a Mahlo cardinal if the
  set of regular cardinals below it is stationary.}), Banach \cite{Banach1930} and Ulam
\cite{Ulam1930} (measurable cardinals\footnote{$\kappa$ is measurable
  if it has a a non-principle ultrafilter}.), Erd\H{o}s and Tarski
\cite{ErdosTarski1943}, \cite{ErdosTarski1961}
(weakly compact cardinals). 
In 1964-1964,  H.~Jerome Keisler and Alfred Tarski
\cite{KeislerTarski1963} 
made a systematic study of weakly compact, measurable and 
strongly compact cardinals.
Supercompact cardinals\footnote{A cardinal $\kappa$ is supercompact
if and only for every $\lambda\ge\kappa$, there is a normal fine ultrafilter $U$
on $\mathcal{P}_\kappa(\lambda)$. Alternatively, a cardinal $\kappa$
is $\gamma$-supercompact for $\gamma\ge\kappa$, if and only if there
is an elementary embedding $j:V\to M$ such that $\kappa$ is the
critical point of $j$, $\gamma<j(\kappa)$, and the model $M$ contains
all of its $\gamma$-size subsets, and is supercompact if and only if
it is $\gamma$-supercompact for all $\gamma\ge\kappa$.} were
introduced by Solovay and William Reinhardt
\cite{ReinhardtSolovayKanamori1978} no later than 1966-67 (see
\cite[page 186]{Magidor1971a}). 
Another strand of large cardinal properties came out of generalization
of partition properties, e.g. Ramsey cardinals\footnote{A cardinal
  $\kappa$ is \emph{Ramsey} if for every coloring by $f$ of its finite
  subsets with two colors, 
there is a subset $H\subseteq\kappa$ of cardinality $\kappa$ 
such that for each positive $n<\omega$, 
all the $n$-element subsets of $H$ receive the same color.}
introduced in 1962
by Erd\H{o}s and Hajnal \cite{ErdosHajnal1962Ramsey}.  

Models are foundational for set theory, since
in forcing one starts with a model and extends it to get a 
new model, as Cohen extended G\"odel's 
Constructible Universe to a model in which the Continuum Hypothesis fails. 
Alfred Tarski and his students developed model theory in the 1950's and 1960's
and his students, C.C.~Chang and Keisler \cite{ChangKeisler1973} 
envisioned in 1963 their textbook on model theory, which would be 
based on lecture notes, with a significant revision after
the 1967 UCLA set theory meeting, but did not appear until 1973. 
Model theoretic techniques were applied in varied ways
to set theory, including the use of absoluteness to transfer results
from one model to another as done by Jack Silver\footnote{Silver received his doctorate from the University of California 
Berkeley in 1966, and joined the faculty there shortly after with
Karel Prikry graduating in 1968 as his first student.}  in his 
1966 thesis \cite{Silver1966PhDthesis}. 

Since forcing employed partial orders (e.g. Cohen reals), Boolean
algebras (e.g. random reals), and trees (e.g. Sacks forcing), they
also became combinatorial objects of study in addition to the graphs and
hypergraphs of the partition calculus.   Notions of largeness included
closed unbounded subsets, stationary subsets. 

\bigskip

Baumgartner was quickly brought up to speed on current topics in
set theory at the four week long 1967 UCLA Summer Institute on Axiomatic Set
Theory. Scott\footnote{Dana Scott presented the Boolean approach
  and was expected to submit a paper on it with Solovay to the proceedings
  of the conference (see \cite{Shoenfield1971}) but did not do so.} 
and Joseph Shoenfield \cite{Shoenfield1971}
gave ten lectures each on forcing. Sacks spoke on the perfect set
forcing or tree forcing named for him.
Many of the other topics that came up at the 
meeting and in the two volume proceedings
\cite{UCLAvol1}, \cite{UCLAvol2}
are related to Baumgartner's published work.
A variety of large cardinals were discussed including 
measurable cardinals, real-valued measurable
cardinals, and supercompact cardinals, as well as reflection principles.
Other topics included $\lambda$-saturated ideals,
extensions of Lebesgue measure, 
the partition calculus, Kurepa's Hypothesis, and Chang's Conjecture.

\medskip

In \cite[161]{MooreG1988} Gregory Moore, based on an
interview with Baumgartner in 1980, reported that ``At Berkeley a
group of young graduate students (including Baumgartner, Laver, and
Mitchell) organized their own seminar --- with no faculty invited.''

In the acknowledgments section 
of his thesis, Baumgartner \cite[2]{JEB1970thesis} 
Baumgartner asserted that his 
``greatest mathematical debt is to the work of Paul Cohen,
and to the work of Robert Solovay and others in making it understandable.''
Baumgartner credited Jack Silver with his 
``initiation into the techniques of forcing proofs,''
and noted that ``most
of the problems treated 
here were suggested to me by Fred Galvin,\footnote{Galvin held pre- and 
post-doctoral positions at University of California Berkeley during 1965-1968} 
Richard Laver,\footnote{Laver \cite{Laverthesis1969} earned his doctorate from 
University of California Berkeley in 1969.} William
Mitchell,\footnote{Mitchell \cite{Mitchell1970} 
received his doctorate in 1970.} and Jack Silver,
and conversations with them have resulted in
the improvement of many proofs and the extension of many results.''  All of those mentioned attended the 1967 Summer School.

\subsection{On Suslin's Question}
In 1970, Baumgartner, Jerome Malitz, and William Reinhardt 
showed, that if the usual axioms of set theory are consistent,  then
so is a positive answer to Mikhail Suslin's
Question of 1920 \cite{Suslin1920}, rephrased below in modern terminology:
\begin{quote}
Must every complete, dense in itself, linear order 
without endpoints for which every pairwise disjoint set of intervals is 
countable be a copy of the real line?
\end{quote}
Baumgartner, Malitz and Reinhardt built on work by 
{\Dbar}uro Kurepa%
\footnote{The masterful survey trees and linear orders
by Stevo Todorcevic \cite{Todor1984settop} 
includes an excellent introduction to Kurepa's work.} 
\cite{Kurepa1935} who 
conducted the first systematic investigation of uncountable trees,
introducing the partition tree of a linear order, the linearization
of a tree, Suslin, Aronszajn, and Kurepa trees%
\footnote{These are all trees of cardinality $\omega_1$ and height $\omega_1$:
a  \emph{Suslin tree} has no uncountable branch and no uncountable 
antichain, where an \emph{antichain} in a partial order is a set whose 
elements are pairwise incomparable; an \emph{Aronszajn tree} has no
uncountable chain and countable levels; and a \emph{Kurepa tree} has countable 
levels and  more than $\aleph_1$ branches.}
He showed the equivalence of the existence of a 
Suslin tree to a negative answer to Suslin's question. 

\smallskip

In their 1970 paper, Baumgartner, Malitz and Reinhardt%
\footnote{Baumgartner (using Martin's Axiom in his thesis) 
and the team of Jerome Malitz and William Reinhardt independently
proved these results.}
proved the existence of a forcing extension in which
$2^{\aleph_0}=\aleph_2$ and 
all Aronszajn trees are embeddable in the rationals,%
\footnote{Kurepa \cite{Kurepa1937Atree} constructed the first Aronszajn tree with
an embedding into the rational numbers.}
where a tree $(T,<_T)$ \emph{embeds in the rationals}
if there is a function $f:T\to\mathbb{Q}$ such that $s<_T t$ implies 
$f(s)<f(t)$.  In such a case we say $(T,<_T)$ is \emph{special}. 

The heart of the argument is
an elegant proof of the countable chain condition of the forcing.

To see how this result is connected to Suslin's Question, note that
if $(T,<)$ is Aronszajn tree with an embedding $f:T\to\mathbb{Q}$,
then $(T,<)$ is 
a union of countably many antichains,%
\footnote{Both being a countable
  union of antichains and having an embedding into $\mathbb{Q}$ have
  been used as the definition of special. See \cite[Exercise
    9.9]{Jech2003} for the equivalence.}
since for each rational $r$, 
$f^{-1}\{ r\}$ is an antichain, 
and it follows that $(T,<_T)$ has an uncountable antichain, since $T$
is uncountable. Thus if an Aronszajn tree is special, it fails to be 
a Suslin tree. By definition, Suslin trees are Aronszajn trees with 
no uncountable antichain, so in the 
Baumgartner-Malitz-Reinhard extension, there are no Suslin trees
giving the consistency relative to ZFC of a positive answer to 
Suslin's Question.  
At the end of the three author paper there were 
three questions, the last of which was asked by Baumgartner: 
Is it consistent with ZFC + $2^{\aleph_0}=\aleph_1$ to 
assume that every Aronszajn tree is embeddable in the rationals?

Baumgartner, Malitz, Reinhardt were preceded by
Robert Solovay and Tennenbaum \cite{SolovayTennenbaum1971} who
found their proof in June 1965 that a positive answer to Suslin's 
Question was relatively consistent with the usual axioms of set theory.  Independently,
Thomas Jech \cite{Jech1967Suslin} and Stanley Tennenbaum 
\cite{Tennenbaum1968Suslin} 
used forcing to show the consistency relative to the usual axioms of
set theory of a negative answer 
to the famous question by Suslin, so a positive answer is independent
of the usual axioms of set theory.  

\medskip

\subsection{Generalized Ramsey Theory}
Next we turn Baumgartner's work in the partition calculus which grew out 
of generalizations of Ramsey's Theorem of 1930. 
Frank P. Ramsey \cite{Ramsey1930} proved that 
for any partition of the $n$-element subsets of an infinite set $A$,
there is an infinite subset $H\subseteq A$, all of whose $n$-element
sets lie in the same cell of the partition.
To present this in modern notation, we introduce the arrow notation 
of Richard Rado \cite{ErdosRado1953}:  for any cardinal $\kappa$, for ordinals
$\langle\alpha_i\mid i<\kappa$ and $\beta$, and any $r\in\omega$  the \emph{partition property} 
\[\beta\rightarrow (\alpha_i)^r_\kappa\]
is the statement that for any partition $f:[\beta]^r\to\kappa$, there
is an $i<\kappa$ and a
subset $A\subseteq \beta$ of order type $\alpha_i$ (in symbols, 
$\otp(A)=\alpha_i$) \emph{homogeneous} for the partition, that is, 
all $r$-tuples from $A$ are in the same cell, i.e. $f$ is constant on $[A]^r$.
In such a case we often call $f$  a \emph{coloring}, $\beta$ a
 \emph{resource} and each $\alpha_i$ a \emph{goal}.  If all the
 $\alpha_i$s are equal, say to $\alpha$, we abbreviate the notation
to $\beta\rightarrow (\alpha)^r_\kappa$.
With this notation in hand, Ramsey's Theorem is the statement that for
all $k<\omega$, $\omega\rightarrow(\omega)^n_k$.  In 1930
Sierpi\'nski proved $\omega_1\nrightarrow (\omega_1)^2_2$.
One of the many equivalent definitions of ``$\kappa$ is a \emph{weakly compact 
cardinal}'' is that $\kappa\rightarrow(\kappa)^2_2$.

In 1973, Baumgartner and Andr\'as Hajnal \cite{JEB1973Hajnal} solved the
$\rho=0$ case of Problem 10 and all of Problem 10A  
of the paper by Erd\H{o}s and Hajnal \cite{EH71.prob}\footnote{This
  paper was based on the lecture by Erd\H{o}s at the 1967 Summer School at UCLA.}
by proving that for all countable ordinals $\alpha$ and finite $k$, 
$\omega_1\rightarrow (\alpha)^2_k$.

According to Hajnal,  
work on what is known as the Baumgartner-Hajnal
Theorem \cite{JEB1973Hajnal} started  in 1970.
Hajnal \cite{Hajnal2013} learned about Martin's Axiom from 
Istv\'an Juh\'asz in Budapest late in 1970.  He decided to try it out 
on the Erd\H{o}s problem $\omega_1\rightarrow(\alpha,\alpha)^2$ and 
was delighted to discover it worked.%
\footnote{Hajnal was not the first to use Martin's Axiom this way.
In a personal conversation, Laver told me that he 
used  Martin's Axiom to prove a partition relation equivalent to
$\omega_1\rightarrow (\omega_1, (\omega:\omega_1))^2$ early in his
post doc at Bristol 1969-1971.
} 
  
Shortly thereafter he attended the International
Congress of Mathematicians in Nice September 1-10, 1970, where he went
around telling people, including Solovay of his proof, but there was
little interest in the result.  Then he contacted Fred Galvin in Budapest and
Galvin suggested getting in touch with Baumgartner.  Later Galvin
wrote to Hajnal that Baumgartner said that he could prove the theorem
outright because there was an argument in Silver's thesis
\cite{Silver1966PhDthesis} (see also \cite{Silver1970} cited in their paper)
that could be used 
to eliminate any appeal to Martin's Axiom by absoluteness.

There was a meeting organized by the New York Academy
of Sciences, and Erd\H{o}s, Hajnal and Baumgartner met there
and talked about the result.  Each told the other what they knew. It
took quite awhile for them to be convinced that both parts were right.
Then Baumgartner wrote it up.  They wrote an initial technical report
\cite{JEB1973Hajnal} published in April 1971 and submitted their final
report the same year but it did not appear until 1973.

To set the Baumgartner-Hajnal Theorem for $\omega_1$ in context, 
note that in 1933, Sierpi\'nski \cite{Sierpinski1933} proved the
analog of Ramsey's Theorem fails for $\omega_1$:
$\omega_1\nrightarrow (\omega_1)^2_2$.  In 1942, Erd\H{o}s
\cite{Erdos42} proved an early positive result in the partition 
calculus for uncountable cardinals when he proved a graph theoretic
equivalence of $\kappa^{\kappa}{}^+\rightarrow (\kappa^+)^2_\kappa$;
a similar theorem was implicit in work of Kurepa \cite{Kurepa1939V} from 1939.

In 1956, Erd\H{o}s and Rado \cite{EP1956Rado} published
the  first systematic treatment of the partition calculus. 
They proved that for any finite $n$, 
$\omega_1\rightarrow (\omega+n)^2_2$ and 
$\omega_1\rightarrow (\omega+1, \omega_1)^2_2$.  They further proved
that for any finite positive $n$, $m$, and $k$, and
any uncountable order type $\varphi$ for which neither
$\omega_1$ nor its reverse, $\omega_1{}^\ast$, embeds in $\varphi$,
 the following partition 
relations hold:
\[ 
\varphi\rightarrow (\omega+n,\omega\cdot n)^2, \quad
\varphi\rightarrow (\omega+n)^2_3, \quad
\varphi\rightarrow (\omega+1)^2_k.
\]
In 1960, Hajnal proved that for $n<\omega$,  
$\omega_1\rightarrow(\omega\cdot 2, \omega\cdot n)$ 
and that under CH, $\omega_1\nrightarrow (\omega+2,\omega_1)^2$.
He further proved that for any uncountable order type $\varphi$ 
for which neither
$\omega_1$ nor its reverse, $\omega_1{}^\ast$, embeds in $\varphi$, 
and any finite positive $n$, countable $\alpha$, and for $\eta$
the order type of the rationals, the following partition 
relations hold:
\[ 
\varphi\rightarrow (\alpha\lor\alpha^\ast,\eta)^2, \quad \mbox{ and }
\varphi\rightarrow (\omega\cdot n,\alpha)^2,
\]
where the goal written $\alpha\lor\alpha^\ast$ is met if
there is a homogeneous set for that color isomorphic to one of
$\alpha$ and $\alpha^\ast$.

Galvin (unpublished) proved no later than May 1970%
\footnote{See \cite{EH74.prob} for the timing.} 
that $\varphi\nrightarrow (\omega)^1_\omega$ implies
$\varphi\nrightarrow (\omega,\omega + 1)^2$. Galvin then  
revised the  conjectures by Erd\H{o}s and Rado in Problems 10, 10A, 11
for $\omega_1$, $\lambda$, the order type of the set of real numbers,
and order types which embed neither $\omega_1$ nor its reverse $\omega_1{}^\ast$
to a conjecture for  order types $\varphi$
for which $\varphi\rightarrow (\omega)^1_\omega$, i.e. 
order types with the property that for every
partition into countably many sets, there is one which includes an
increasing sequence.   

The full Baumgartner-Hajnal Theorem asserts the revised conjecture is
true: for any order type $\varphi$, if $\varphi\rightarrow
(\omega)^1_\omega$, then for all $\alpha<\omega_1$ and $k<\omega$, 
\[\varphi\rightarrow (\alpha)^2_k.\]

Its metamathematical proof took a result proved with additional
assumptions and then showed the result is absolute, that is, its truth in a
model  ZFC with additional assumptions implies 
that it is a consequence of ZFC. This paper introduced this method
to a wide audience.

One of the key lemmas in the proof is the preservation result which
says that if $\varphi$ is an order type such that
$\varphi\rightarrow(\omega)^1_\omega$, then $V^\mathbb{P}\models
\varphi\rightarrow(\omega)^1_\omega$ for every ccc forcing notion $\mathbb{P}$.

Subsequently Galvin \cite{Galvin1975BH} gave a 
proof of the Baumgartner-Hajnal Theorem that was purely combinatorial. 
Confirming a conjecture of Galvin, Todorcevic \cite{Todor1985poset} extended the
Baumgartner-Hajnal Theorem to the class of all partially ordered sets by
proving that for every partial order $P$, if $P\rightarrow (\omega)^1_\omega$,
then $P\rightarrow (\alpha)^2_k$ for all $\alpha<\omega_1$ and $k<\omega$.
Note that in this context the analogous absoluteness result that
$P\rightarrow(\omega)^1_\omega$ is preserved by ccc forcing is no
longer true, so \cite{Todor1985poset} used a different argument.
In 1983, Todorcevic \cite{Todor1983forcePos} through forcing showed the 
consistency of a partition relation considerably stronger than the 
Baumgartner-Hajnal Theorem:\footnote{Todorcevic
  \cite{Todor1983forcePos} reported the work on the paper was done
  during the academic year 1980-1981 when he was visiting Dartmouth
  College.} 
$\omega_1\rightarrow (\omega_1,\alpha)^2$ for all countable ordinals $\alpha$.

In 1991, Prikry and  Milner \cite{MilnerPrikry1991}
proved $\omega_1\rightarrow (\omega\cdot 2 +1, 4)^3$
by first showing  that the partition relation holds in a model of
Todorcevic in which both Martin’s Axiom and 
$\omega_1\rightarrow (\omega\cdot 2 +1)^2$ hold, and then using the
approach taken by Baumgartner and Hajnal to show the consistency
result is actually a ZFC theorem.

\subsection{Basis Problem for Uncountable Order Types}
A set of order types $B$ is a \emph{basis} for a family 
$\mathcal{F}$ of linear orders if $B\subseteq\mathcal{F}$ and for 
all $\varphi\in\mathcal{F}$ there is some $\psi\in B$ with $\psi$
embeddable in $\varphi$, in symbols $\psi\le\varphi$.
For example, $B_{\aleph_0}=\{\omega, \omega^*\}$ is a two-element 
basis for countably infinite linear orders,  and $B_{\mathbb{Q}}=\{\eta\}$ is a
one-element basis for all infinite linear orders that are dense in themselves, since $\eta$, the order type of the rationals
is embeddable in every infinite dense in itself order type.  

A set $A$ of real numbers is $\aleph_1$-dense if it has cardinality
$\aleph_1$ and between any two
elements of $A$ there are exactly $\aleph_1$ members of $A$.
In 1973
Baumgartner \cite{JEB1971CantorThm.abs}, \cite{JEB1973reals}  
published his best known theorem on order:
it is relatively consistent with ZFC that
  $2^{\aleph_0}=\aleph_2$ and all $\aleph_1$-dense sets of reals are
  isomorphic.
A key idea of the proof is to start with a model of the
Continuum Hypothesis, build an iteration of 
countable chain condition  (ccc) forcings 
that introduces the necessary order isomorphisms between pairs of 
$\aleph_1$-dense  sets. The CH is preserved at each intermediate stage 
and is used to prove that the forcings adding
isomorphisms satisfy the ccc.  At the end of the iteration, the
Continuum Hypothesis no longer holds, but its presence in the
intermediate stages was sufficient.  

Another approach to Baumgartner's result (using the PFA context) is given in
Chapter 8 of Todorcevic's book \cite{TodorProbTop1989}, and  
yet another approach to Baumgartner's result are generalizations given
in \cite{AbrahamRubinShelah1985}.

In 1980, Shelah \cite{Shelah1980LinO}
proved the consistent existence of a linear order universal in power $\aleph_1$ with the continuum equal to $\mathfrak{c}=\aleph_2$, 
Shelah compared and contrasted his use of oracle forcing in the proof with
that of Baumgartner in his proof that all $\aleph_1$-dense sets of reals
can be isomorphic.  

Baumgartner used his forcing result to show the consistency of ZFC + ``the class of real 
types has a one element basis.'' He also asked if it is consistent for
all $\aleph_2$-dense sets of reals to be isomorphic,
and the question was answered in the positive by Itay Neeman (email of 
March 5, 2016). 

Baumgartner noted his proof can be extended to add Martin's Axiom to 
the conclusion,
and asked if ``all $\aleph_1$-dense sets of reals are isomorphic'' follows
from Martin's Axiom + $2^{\aleph_0}>\aleph_1$.  
Uri Abraham and Saharon Shelah 
gave a negative answer in 1981, and Abraham, Matatyahu Rubin, and Shelah 
\cite{AbrahamRubinShelah1985}
showed that it is relatively consistent with ZFC that
$2^{\aleph_0}>\aleph_2$ and all $\aleph_1$-dense sets of reals are isomorphic
and proved that if all $\aleph_1$-dense sets of reals are isomorphic, then
$2^{\aleph_0}<2^{\aleph_1}$.   

\medskip

In 1976 Baumgartner  \cite{JEB1976newotp} considered the problem of finding 
a nice basis for the class $\Phi$ of all uncountable order types which cannot 
be represented as the union of countably many well-orderings.  He set 
$\Phi_1=\{\omega_1{}^\ast\}$, let 
$\Phi_2$ be the uncountable order types embeddable in the reals, and let
$\Phi_3$ be the uncountable order types which do not embed a subset of
type $\omega_1$ nor a subset of type $\omega_1{}^*$ nor an uncountable
subset of the real numbers. He called the type types in $\Phi_3$ \emph{Specker types}.\footnote{In 
\cite[page 443]{ErdosRado1956} Erd\H{o}s and Rado  conjectured that there were no 
uncountable linear order types which did not embed $\omega_1$, nor embed 
the reverse, $\omega_1{}^\ast$, nor embed an uncountable subset of the real 
numbers, but they included a footnote that Specker had disproved the 
conjecture.} 
Finally he let $\Phi_4$ be the uncountable order types $\varphi$ such that 
every uncountable subtype $\psi\le\varphi$ contains an uncountable 
well-ordering but $\varphi$ cannot be represented as the union of 
countably many well-orders.  Then every element of $\Phi$ embeds some
element of $\Phi_1\cup\Phi_2\cup\Phi_3\cup\Phi_4$,
and the question of a basis for $\Phi$ can be subdivided into finding
bases for the components. Consequently, a basis for uncountable order 
types can be obtained by adding $\omega_1$ and its reverse, $\omega_1{}^\ast$ 
to a basis for $\Phi_2\cup \Phi_3$, since every element of $\Phi_4$ embeds 
$\omega_1$.  Baumgartner answered Galvin's question of whether
$\Phi_1\cup\Phi_2\cup\Phi_3$  formed a basis for $\Phi$ 
by showing that $\Phi_4$ is non-empty.  

Most of the paper is devoted to developing a structure theory for
elements of $\Phi_4$, where stationary sets play a significant role.
In particular, in Corollary 7.9, he proved that 
if ZFC + ``there exists a weakly compact cardinal'' 
is consistent, then so is ZFC + 
``for every stationary subset $C\subseteq\omega_2$, if
$\operatorname{cf}(\alpha)=\omega$ for all $\alpha\in C$, then there
is an $\alpha<\omega_2$ such that $C\cap\alpha$ is stationary in
$\alpha$.''    
In modern language, the latter statement translates to every stationary subset 
of $\omega_2\cap \cof (\omega)$ reflects,\footnote{%
In 1985 Leo Harrington and Shelah \cite{HarringtonShelah1985} showed
that ZFC + ``the  existence of a Mahlo cardinal'' is equiconsistent with
ZFC + ``every stationary subset of $\omega_2\cap\cof(\omega)$ $\omega$ reflects.''}
 where $\cof (\omega)$ is the
collection of ordinals of cofinality $\omega$.
In \cite{JEB1976newotp} Baumgartner (as reworked later by Shelah and
perhaps others) showed that if one collapses a supercompact cardinal
to $\omega_1$ via a standard L\'evy collapse, then in the resulting
forcing extension, for every regular cardinal $\kappa\ge\omega_2$, and
every stationary set $S\subseteq \kappa\cap\cof(\omega)$, there is a
$\gamma<\kappa$ such that $S\cap\gamma$ is stationary in $\gamma$.
Menachem Magidor \cite[page 756]{Magidor1982} asserted that
Baumgartner's proof actually showed under the hypothesis of Corollary 7.9, that ZFC + ``every pair
of stationary subsets of $\omega_2\cap\cof(\omega)$ has a common point
of reflection.''  Magidor also showed the equiconsistency of this
statement with the hypothesis of Corollary 7.9. 

In his paper on uncountable order types,
Baumgartner also assembles the ingredients for a proof that 
every element of $\Phi_3$ is the linearization of an Aronszajn tree, 
that is an \emph{Aronszajn line}, so $\Phi_3$ and the collection of 
Aronszajn lines coincide.\footnote{Todorcevic 
\cite{Todor1984settop} described this result as part of the folklore of
the subject, noting that a large part of it was proved by Kurepa
\cite[page 127-9]{Kurepa1935}, and further referring the interested
reader to the Erd\H{o}s-Rado paper \cite{ErdosRado1956} which only has 
their conjecture and the note that Specker refuted it, and 
to a survey paper by R.~Ricabarra.
Proofs of this fact can be found in Baumgartner's 1982 survey
article \cite{JEB1982order} on order and Todorcevic's 1984
survey article \cite{Todor1984settop} on trees and linear orders.}

At the end of the paper Baumgartner \cite{JEB1976newotp} 
asked in Problem 5(i)  if ZFC + ``$\Phi_3$ (Aronszajn lines) 
has a finite basis'' is consistent 
and in Problem 5(ii) if ZFC + ``$\Phi_2\cup\Phi_3$ (real types and 
Aronszajn lines) has a finite basis'' is consistent.  
We will return this question in the  section on forcing for all.

\medskip

\subsection{Disjoint Refinements}
In 1975, Baumgartner, Hajnal and Attila M\'at\'e \cite{JEB1975HMsaturation} 
gave a partial answer to a question of Fodor by giving a condition on
the non-stationary ideal $\mbox{NS}_{\omega_1}$
which guarantees any $\omega_1$-sequence of 
stationary subsets of $\omega_1$ can be resolved into an antichain in 
$\mathcal{P}(\omega_1)/\mbox{NS}_{\omega_1}$, as described in the following theorem.

\begin{theorem}[Baumgartner, Hajnal, M\'at\'e] Assume $I$ is a normal ideal on 
$\omega_1$ such that given any 
$\langle X_\alpha\in\mathcal{P}(\omega_1)\setminus I\mid\alpha<\omega_1\rangle$,
there exists $X\in\mathcal{P}(X_0)\setminus I$ such that 
$X_\alpha\setminus X\notin I$ for each $\alpha\ge 1$.\footnote{%
This condition holds if $Y_{\omega_1}$ fails to be $\omega_1$-dense,
where $\mathcal{I}$ is \emph{$\omega_1$-dense} on $\omega_1$ if there 
is a family $D\subseteq \mathcal{P}(\omega_1)\setminus \mathcal{I}$ 
of power $\omega_1$ such that for any $X\in\mathcal{P}\setminus I$, 
there is a $Y\in D$ with  $Y\setminus X\in \mathcal{I}$.}  
Then, given
$\langle S_\alpha\in\mathcal{P}(\omega_1)\setminus I\mid \alpha<\omega_1\rangle$,
there exist $\langle A_\alpha\in\mathcal{P}(\omega_1)\setminus I\mid\alpha<\omega_1\rangle$ 
such that each $A_\alpha\subseteq S_\alpha$ and the $A_\alpha$'s are disjoint.
\end{theorem}
In 2000, Paul Larson \cite{Larson2000} used the above theorem   
in his article \emph{Separating Stationary Reflection Principles}
to show that  Todorcevic's Strong Reflection Property (SRP) implies 
SR${}_{\omega_1}^\ast$, the strongest strengthening of Stationary Reflection (SR)
that Larson considered.

In 2014, Monroe Eskew \cite{Eskew2014PhDthesis} showed that the
Baumgartner-Hajnal-M\'at\'e result cannot be lifted to ideals in general on 
larger $\aleph_n$s 
when he proved 
that if ZFC + ``there is an almost huge cardinal'' is consistent,
then for $n>1$, 
so is ZFC + GCH + ``there is a normal $\aleph_n$-complete 
ideal $\mathcal{I}$ on $\aleph_n$ 
and a sequence of $\aleph_n$ 
many $\mathcal{I}$-positive sets  which has no disjoint refinement''
is consistent. 

\subsection{Almost Disjoint Families}
In 1976 Baumgartner \cite{JEB1976ad} extended the work \cite{Sierpinski1928ad}, \cite{Tarski1928ad}, 
\cite{Tarski1929ad} of Sierpi\'nski and Tarski\footnote{They assumed GCH.}
to give a complete solution to their questions under the assumption
of the Generalized Continuum Hypothesis (GCH) 
after connecting questions about families of almost disjoint sets 
with his \emph{dense set problem} 
of finding out, for a given
cardinal $\kappa$, which cardinals $\lambda$ support a linear order
with a dense subset of power $\kappa$.%
\footnote{Baumgartner formulated the dense set problem
by generalizing a question of Malitz 
\cite{Malitz1968} in a proof that the Hanf number for 
complete $L_{\omega_1,\omega}$ sentences is $\beth_{\omega_1}$.
Baumgartner \cite{JEB1974Hanf} 
eliminated the use by Malitz of GCH by applying a combinatorial fact due to Hausdorff.}
William Mitchell \cite{Mitchell1972Atree} 
also worked on the dense set problem, and Baumgartner documented the 
interconnections of their results.  
Baumgartner developed extensions of an observation of Shelah 
that permit transfer of a result about the existence of a set $S$ of power
$\lambda$ with a dense subset $U$ of power $\kappa$ to the statement
obtained by replacing $\lambda$ by $\lambda^\rho$ and $\kappa$ by $\kappa^\rho$,
and drew conclusions from the extensions about the existence of various
families of almost disjoint sets. He constructed
a broad range of forcing extensions modeling a variety of answers
to the almost disjoint sets questions and the dense sets problems
using along the way the Erd\H{o}s-Rado Theorem for partitions with many parts,
Jensen's $\lozenge_\kappa$, and Easton forcing.
In particular, he proved that it is consistent with ZFC that 
$2^{\aleph_0}=\aleph_{\omega_1}$, 
$2^{\aleph_1}=\aleph_{\omega_1+1}$ and there is no family of $2^{\aleph_1}$ pairwise
almost-disjoint subsets of $\aleph_1$.   

As a special case of a more general theorem, 
Sierpi{\'n}ski \cite{Sierpinski1928ad}, using the Generalized Continuum 
Hypothesis, proved that for any infinite set $A$, there is a family of 
size $\aleph_1$ of strongly almost disjoint subsets of $A$. 
In 1934, Sierpi{\'n}ski  labeled this proposition $P_{11}$,
and showed it is equivalent to the Continuum Hypothesis in his
book \cite{Sierpinski1956}.  Baumgartner \cite[pages 424, 428]{JEB1976ad} proved 
consistency with ZFC and independence from ZFC of the existence of 
a strongly almost disjoint family of uncountable subsets of $\omega_1$.
He started with an almost disjoint family $\mathcal{F}$ of $\aleph_2$ many 
uncountable subsets of $\omega_1$,\footnote{Such a family can be constructed 
by transfinite recursion.} enumerated the family as $\langle F_\alpha\mid \alpha
<\omega_2$, and let $\mathbb{P}$ be the set of all finite partial functions 
$f: \omega_2 \to[\omega_1]^\omega$ such that for all $\alpha\in\dom (f)$, $f(x)$ 
is a finite subset of $F_\alpha$.  He let $f\le g$ if and only if (a) $\dom(g)
\subseteq \dom(f)$; (b) for all $\alpha$ in the domain of $g$, 
$g(\alpha)\subseteq f(\alpha)$, and (c) for all $\alpha\ne\beta$ in the domain 
of  $g$, $f(\alpha)\cap f(\beta)=g(\alpha)\cap g(\beta)$.  In the forcing
extension $\mathfrak{M}[G]$, the sets $G_\alpha=\bigcup\{f(\alpha)\mid f\in G\}$
are considered to be obtained by \emph{thinning out} the $F_\alpha$, and 
they form the strongly almost disjoint family in the extension. 

Baumgartner used results and adapted  techniques used in the study of 
almost disjoint families to prove theorems about polarized partitions 
in the final section of his paper.
For cardinals $\kappa$, $\lambda$, $\mu$, $\nu$, and $\rho$,
the \emph{polarized partition relation} 
{\small $\left(\begin{array}{c}\kappa\\ \lambda\end{array}\right)\rightarrow 
 \left(\begin{array}{c}\mu\\ \nu\end{array}\right)_\rho$}
holds if and only  if for all $f:\kappa\times\lambda\rightarrow\rho$ 
there are 
$A\subseteq \kappa$ with $\otp(A)=\mu$ and 
$B\subseteq \lambda$ with $\otp (B)=\nu$
such that $f$ is constant on $A\times B$.  We also consider the variant
where the subscript is $<\rho$ for coloring maps whose range has
cardinality $<\rho$, and the variant
{\small $\left(\begin{array}{c}\kappa\\ \lambda\end{array}\right)\rightarrow 
 \left(\begin{array}{cc}\mu&\sigma\\ \nu&\tau\end{array}\right)$}
when there are only two color classes and different goals for the
different colors.
  This relation was introduced in \cite{EP1956Rado} and studied in 
\cite{EP1965HajnalRado}.  Baumgartner results include the following
where $\mathfrak{c}=2^{\aleph_0}$:
\begin{enumerate}
\item
{\small $\left(\begin{array}{c}\mathfrak{c}\\ \aleph_1\end{array}\right)\rightarrow 
 \left(\begin{array}{c}\mathfrak{c}\\ \alpha\end{array}\right)$} for
all $\alpha<\omega_1$, but {\small $\left(\begin{array}{c}\mathfrak{c}\\ \aleph_1\end{array}\right)\nrightarrow 
 \left(\begin{array}{c}2\\ \aleph_0\end{array}\right)_\omega$}.
\item
{\small $\left(\begin{array}{c}\mathfrak{c}\\ \aleph_1\end{array}\right)\rightarrow 
 \left(\begin{array}{cc}\mathfrak{c}&\mathfrak{c}\\ \aleph_1&\alpha\end{array}\right)$} for
all $\alpha<\omega_1$, but 
{\small $\left(\begin{array}{c}\mathfrak{c}\\ \aleph_1\end{array}\right)\nrightarrow 
 \left(\begin{array}{cc}\aleph_1&\aleph_1\\ \aleph_0&\aleph_0\end{array}\right)_\omega$}.
\end{enumerate}
 
Almost disjoint families have been and continue to be used to construct 
interesting examples.  Haim Gaifman and Specker \cite{GaifmanSpecker1964} 
showed in 1964 that if $\kappa^{<\kappa}=\kappa$, then there are $2^{\kappa^+}$ many
different types of normal $\kappa^+$-Aronszajn trees by using
a family of almost disjoint sets in their construction.

In 2005, using results under GCH of Sierpi\'nski and Tarski cited above,
Lorenz Halbeisen \cite{Halbeisen2005} proved 
the consistency with ZFC that for all cardinals $\kappa$, every infinite
dimensional Banach space of cardinality $\kappa$ admits $2^\kappa$ pairwise
almost disjoint normalized Hamel bases. By way of contrast, using a result 
of Baumgartner \cite[Theorem 5.6(b)]{JEB1976ad}, Halbeisen proved the 
consistency with ZFC that $2^\kappa\le \kappa^{++}$ and no infinite dimensional 
Banach space of cardinality $\kappa$ admits $\kappa^{++}$ pairwise almost 
disjoint normalized Hamel bases.  

In 2006, J. Donald Monk   \cite{Monk2006} revisited 
and extended Baumgartner's work on families of almost disjoint sets with a 
focus on the sizes of maximal families.

Cristina Brech and Piotr Koszmider 
\cite{BrechKoszmider} used the product of Baumgartner's $\mathbb{P}$
with the standard $\sigma$-closed and $\omega_2$-cc forcing for adding
$\omega_3$ subsets of $\omega_1$ with countable conditions
in their construction of a forcing extension in which 
there is no universal Banach space of density the continuum.  

\subsection{Translating stationary to closed unbounded}
In 1976, Baumgartner, Harrington and Eugene Kleinberg \cite{JEB1976HK}
were able to  add a closed unbounded set as a subset of a stationary set
$A\subseteq\omega_1$ by forcing with closed countable subsets of 
$A$ whose order type is a successor ordinal.  This technique is useful
for translating problems about stationary sets into ones about closed
unbounded sets, and closed unbounded sets provide the ladder for
recursive constructions and inductive proofs.

To start, Baumgartner, Harrington and Kleinberg recalled the well-known theorem
that for any regular uncountable cardinal $\kappa$, the intersection
of fewer than $\kappa$ many closed unbounded sets  is 
closed unbounded (we will abbreviate ``closed unbounded'' to \emph{club}.)
Then they observed that $\mathcal{F}_\kappa$, the family of all subsets of
$\kappa$ that have a club set as a subset, is a $\kappa$-additive
non-principal filter on $\kappa$.  

Next they considered the
possibility that $\mathcal{F}_\kappa$ is an ultrafilter. 
If $\kappa>\omega_1$, then $\mathcal{F}_\kappa$ cannot be
an ultrafilter.\footnote{For regular $\kappa>\omega_1$, consider the set 
$A$ of ordinals $\alpha<\kappa$ of cofinality $\omega_1$.  
Neither it nor its complement can contain a closed unbounded 
subset of $\kappa$.}
If $\kappa=\aleph_1$ and $\mathcal{F}$ is an ultrafilter, 
then $\aleph_1$ is a measurable cardinal. Solovay 
\cite{Solovay1971} had shown $\aleph_1$ being measurable
is relatively consistent with the usual axioms of set theory (ZF)
together with the Axiom of Determinacy in his model in which 
all sets of reals are Lebesgue measurable where 
the Axiom of Choice is not true.  Baumgartner, Harrington, and
Kleinberg  then proved that if one
starts with a model $\mathcal{M}$ of the usual axioms of set theory (ZF)
together with the Axiom of Choice and a set $A$ in $\mathcal{M}$ 
such that $\mathcal{M}\models$ ``$A\subseteq\aleph_1$ 
is not disjoint from any closed unbounded subset of $\mathcal{M}$,''
then there is a generic extension $\mathcal{N}$ with the same reals
as $\mathcal{M}$ such that $\mathcal{N}\models$ ``$A$ contains a
closed unbounded subset of $\aleph_1$.''  

Baumgartner, Harrington, and Kleinberg  described their ``shooting a club
through a stationary set'' result as an extension of the theorem of 
Harvey Friedman \cite{FriedmanH1974}, who proved that every stationary
subset of $\aleph_1$ contains arbitrary long countable closed sets.  

In 1978, Abraham and Shelah \cite[page 647-8, Theorem 4]{AbrahamShelah1983}, building on work
by Jonathan Stavi,\footnote{Stavi was cited for handwritten notes from
  1975 on \emph{Adding a closed unbounded set}.} 
proved that if $\kappa=\mu^+$, $\mu^{<\mu}=\mu$, 
and $S$ is a fat stationary 
set,\footnote{A stationary set $S\subseteq \kappa$ is called \emph{fat} 
if and only if for every closed unbounded set $C\subseteq\kappa$, $S\cap C$ 
contains closed sets of ordinals of arbitrarily large order-types below $
\kappa$.  The terminology is from \cite{FleissnerKunen1978}.}
then there is a partial order such that forcing with it introduces a club subset 
of $S$, does not collapse any cardinals, and does not add new subsets of size 
$<\mu$.

In Section 6 of his \emph{Handbook of Set Theory} chapter, 
James Cummings \cite{Cummings2010hand} used this forcing of 
Baumgartner, Harrington and Kleinberg to show that in general 
$(\omega_1,\infty)$-distributivity is weaker than 
$<\omega_1$-strategic closure.

\subsection{Cardinal Arithmetic Constraints}
In 1976 Baumgartner and Karel Prikry \cite{JEB1976Prikry} 
published their elementary proof\footnote{Ronald Jensen (unpublished)
also had an alternative proof.}
of the remarkable result of Jack Silver \cite{Silver1975} 
that if the Generalized Continuum Hypothesis  holds below a 
cardinal $\kappa$ of uncountable cofinality, then it holds at $\kappa$. 
Silver's paper in the proceedings of the International Congress of 
Mathematicians  of 1974 used metamathematical arguments and his four 
page proof omitted many details. Let us note that 
Silver's result has deep roots:  Hilbert 
put the Continuum Problem first on his famous list of problems of
1900.  Felix Hausdorff \cite[page 133]{Hausdorff1907a} speculatively 
used the possible generalization of the continuum hypothesis to larger
$\aleph_\nu$s 
in his analysis of order types generalizing the rationals to larger 
cardinalities. Alfred Tarski \cite[page 10]{Tarski1925gch} 
used the phrase \emph{generalized
continuum hypothesis} (\emph{hypoth\`ese g\'en\'eralis\'ee du continu}) 
in 1925.  
As a service to the wider mathematical community, 
Baumgartner and Prikry \cite{JEB1977Prikry} then wrote an article for
the \emph{American Mathematical Monthly} on the special case 
$2^{\aleph_{\omega_1}}= \aleph_{\omega_1+1}$ 
that used only K\"onig's Theorem (the sum of an indexed family
of cardinals is less that the product of the indexed family), the Regressive 
Function Theorem and basic facts about cardinal arithmetic.

\subsection{Filters, ideals, and partition relations}
Cardinality was the initial notion of largeness for homogeneous sets
in the systematic study of the partition calculus by Erd\H{o}s and 
Rado \cite{EP1956Rado}.  Other notions considered early were
having a large order type, e.g. subsets of the reals order isomorphic
to the set of reals and being in a $\kappa$ complete ultrafilter  
for a large cardinal $\kappa$.

The central point for Baumgartner in his two papers on
ineffable cardinals was that ``many `large cardinal' properties are
better viewed as properties of normal ideals than as properties of
cardinals alone,''  He used a variety of partition relations,
inaccessibility and indescribability in his characterizations of
the normal ideals associated with a variety of mild large
cardinals, that is, ones below a measurable cardinal.  
He also was able to calibrate the strength of the large cardinal
needed for the partition relations in question.
%

In his 1975 paper on ineffable cardinals
Baumgartner \cite{JEB1975.ineff}  analyzed 
large subsets of ineffable, almost ineffable, and subtle cardinals. 
These cardinals  had been introduced by Jensen and Kunen 
\cite{JensenKunen1969} in their analysis of combinatorial principles 
that hold in $L$.

Suppose $A\subseteq\kappa$.
Recall a function $f:A\to\kappa$ is \emph{regressive}
if $f(\alpha)<\alpha$ for all $\alpha>0$;
and $f:[A]^n\to\kappa$ is \emph{regressive} if
$f(\vec a) <\min(\vec a)$ for all $\vec a\in [A]^n$.
In a modern definition, a regular cardinal $\kappa$ is
(1) \emph{ineffable}; (2) \emph{weakly ineffable}, (3) \emph{subtle}
respectively  if and only if for every regressive function 
$f:\kappa\to\mathcal{P}(\kappa)$,
\begin{enumerate}
\item (ineffable)
there is a set $A\subseteq \kappa$ such that the set
$\{\alpha<\kappa\mid A\cap\alpha=f(\alpha)\}$ is stationary;
\item (weakly ineffable)
there is a set $A\subseteq \kappa$ such that the set
$\{\alpha<\kappa\mid A\cap\alpha=f(\alpha)\}$ has cardinality $\kappa$;
\item (subtle) 
for every closed unbounded subset $C\subseteq \kappa$ 
there are $\alpha < \beta\in C$ with $A_\alpha=A_\beta\cap\alpha$.
\end{enumerate}
Kunen \cite{JensenKunen1969} 
proved that a cardinal $\kappa$ is ineffable if and only if
it satisfies the partition relation $\kappa\rightarrow
(\text{stationary})^2_2$, where one asks for a stationary homogeneous
set rather than one of cardinality $\kappa$.  He also proved
that ineffable cardinals were $\Pi^1_2$-indescribable, and 
Kunen and Jensen located the least ineffable cardinal
above the least cardinal which is $\Pi^n_m$-indescribable for all
$m,n<\omega$ and above the least cardinal cardinal $\lambda$ such that
$\lambda\rightarrow (\omega)^{<\omega}_2$.  
They showed weakly ineffable cardinals were 
$\Pi^1_1$-indescribable.  

Baumgartner refined subtle, weakly ineffable, and ineffable to 
$n$-subtle, $n$-weakly ineffable, and $n$-ineffable. 
Extend the notion of regressive to functions whose domain 
is a subset $A\subseteq \kappa$ and whose range is a subset of
$\mathcal{P}(\kappa)$ as follows: 
a function $f:A\to\mathcal{P}(\kappa)$ is \emph{regressive}
if $f(\alpha)\subseteq \alpha$ for all $\alpha\in A$ with $\alpha>0$;
and $f:[A]^n\to\kappa$ is \emph{regressive} if
$f(\vec a) \subseteq \min(\vec a)$ for all $\vec a\in [A]^n$ with
$\min(\vec a)>0$.
Call a set $H\subseteq A$ \emph{set-homogeneous} for a regressive function
$f:[\kappa]^n\to\mathcal{P}(\kappa)$ if and only $H\subseteq A$ and
for all $\vec a, \vec c\in[A]^n$, if $\min(\vec a)\le\min(\vec c)$,
then $f(\vec a)=f(\vec c)\cap\min(\vec a)$. The concept extends in
the natural way for regressive functions from $A$ to $\mathcal{P}(\kappa)$.  

A cardinal  $\kappa$
is $n$-subtle if and only if for every regressive function 
$f:[\kappa]^n\to\mathcal{P}(\kappa)$
and every closed unbounded set $C\subseteq \kappa$, there is a
set $H\in [C]^{n+1}$ homogeneous for $f$.  Also $A\subseteq \kappa$
is $n$-ineffable ($n$-weakly ineffable) if and only if 
every regressive function $f:[A]^n\to\mathcal{P}(\kappa)$ has a homogeneous set
which is stationary in $\kappa$ (of power $\kappa$).  
For any of these cardinal properties, he 
spoke of a subset $A\subseteq \kappa$ as having 
the corresponding homogeneity property, if every suitable regressive
function had the corresponding homogenity property.

Baumgartner proved that for each of the notions of largeness for a cardinal
$\kappa$ the corresponding set of small (i.e. not large) sets forms
a $\kappa$-complete normal ideal on $\kappa$. 

He used partition properties to characterize the various notions 
of largeness. He had equivalences for a subset $A\subseteq\kappa$
being subtle, and being $n$-weakly ineffable using  both
regular and regressive partition relations. For example,
the following are equivalent for a subset $A$
of a regular cardinal $\kappa$:
\begin{enumerate}
\item $A$ is $n$-ineffable.
\item $A\rightarrow (\text{stationary set})^{n+1}_2$.
\item $A\rightarrow (\text{stationary set},\kappa)^{n+2}_2$.
\item $A\rightarrow (\text{stationary set},n+3)^{n+2}_2$.
\end{enumerate}
Thus $\kappa\rightarrow (\text{stationary set})^{m}_2$
does not imply $\kappa\rightarrow (\text{stationary set})^{m+1}_2$,
which stands in contrast to the fact
that $\kappa\rightarrow (\kappa)^2_2$ implies
$\kappa\rightarrow (\kappa)^n_2$ for all positive integers $n$.

Harvey Friedman \cite{FriedmanH2001} has adapted Baumgartner's approach to 
$n$-subtle cardinals for his program to develop ``natural'' propositions 
of finite mathematics whose consistency requires use of large cardinals.
Pierre Matet \cite{Matet2003} has used it to prove a partition property for $\mathcal{P}_\kappa(\lambda)$.

\medskip

In 1977, Baumgartner \cite{JEB1977ineff2}
extended the association of normal ideals 
with large cardinals to weakly compact cardinals and Ramsey cardinals, 
using the notions of
$\alpha$-Erd\H{o}s cardinals%
\footnote{A cardinal $\kappa$ is 
\emph{$\alpha$-Erd\H{o}s} if it is regular and for every regressive function
$f:[\kappa]^{<\omega}\to\kappa$ and every closed unbounded
set $C\subseteq \kappa$, there is $A\subseteq C$ of order type
$\alpha$ which is homogeneous for $f$, i.e. for every positive $n$, 
$f$ is constant on the $n$-element subsets of $A$.}
and canonical sequences.  
Qi Feng \cite{Feng1990} extended 
Baumgartner's work in his thesis, and Ian Sharpe and Philip Welch
\cite{SharpeWelch2011} used Baumgartner's 
canonical sequences and Feng's work to develop an $\alpha$-weakly Erd\H{o}s 
hierarchy as part of their study of implications for inner models
of strengthenings of Chang's conjecture.  They also modeled
the proofs of some of their lemmas on proofs from \cite{JEB1991}.

In 1977 Baumgartner, Taylor and Wagon \cite{JEB1977TW}
used Mahlo's operation $M$ to define a family of ideals 
they called $M$-ideals and to define the notion of a cardinal $\kappa$ being 
\emph{greatly Mahlo}, and then proved that a cardinal was greatly Mahlo if and 
only if it bears an $M$-ideal.  
A very satisfying instance of equiconsistency of 
a combinatorial principle with the existence of a large cardinal was proved
in 2011 by John Krueger and Ernest Schimmerling 
\cite{KruegerSchimmerling2011}. They showed
that the existence of a greatly Mahlo cardinal is equiconsistent with
the existence of a regular uncountable cardinal $\kappa$ such that 
no stationary subset of $\kappa^+$ consisting of ordinals of
cofinality $\kappa$ carries a partial square.\footnote{Partial square
sequences were introduced as a weaking of square sequences
by Shelah (see \cite{Shelah1991partial}).
Suppose $\nu<\kappa^+$ is regular and $A \subseteq
\kappa^+\cap\cof(\nu)$.  Then \emph{$A$ carries a partial square} if
there is a sequence $\langle c_\alpha\mid \alpha\in A\rangle$ such
that (a) each $c_\alpha$ is a closed unbounded subset of $\alpha$ of
order type $\nu$ and whenever $c_\alpha$ and $c_\beta$ share a common
limit point $\gamma$, then $c_\alpha\cap\gamma= c_\beta\cap\gamma$.}

\subsection{Saturated ideals}
Recall that an ideal $\mathcal{I}$ on a cardinal $\kappa$ is
\emph{$\lambda$-saturated}\footnote{Tarski \cite{Tarski1945} introduced  $\lambda$-saturation of ideals in 1945.}
if and only if every pairwise $\mathcal{I}$-almost 
disjoint collection $F\subseteq \mathcal{I}^+$ of $\mathcal{I}$-positive
sets is of cardinality less than $\lambda$. 
Solovay had shown that if a regular cardinal $\kappa$ has a nontrivial normal
$\kappa$-complete $\lambda$-saturated ideal
for some $\lambda<\kappa$, then $\kappa$ is measurable in an inner model. 
In 1970, Kunen \cite{Kunen1970} extended Solovay's result by showing the result was true for regular cardinals $\kappa$ carrying a non-trivial $\kappa$-complete 
$\kappa^+$-saturated ideal.  Kunen noted that it was unknown whether $\omega_1$ can bear an $\omega_2$-saturated ideal. He showed that if a successor cardinal $\kappa$ has a non-trivial $\kappa$-complete $\kappa^+$-saturated ideal then
$0^{\dagger}$ exists.

\medskip

In 1972, Kunen \cite{Kunen1978} proved
that if  ZFC + ``there exists a huge cardinal'' is consistent, then
so is ZFC + ``there is an $\omega_2$-saturated ideal on $\omega_1$.''
In his review of the history of the problem, he noted
that for an uncountable cardinal $\kappa$, the larger the $\lambda$,
the weaker the property of being $\lambda$-saturated ideal on
$\kappa$. The existence of an $\omega$-saturated ideal on
$\kappa$ was equivalent to $\kappa$ being measurable;
the existence of a $\kappa^+$ saturated ideal on $\kappa$
implied $\kappa$ is measurable in an inner model, and the
existence of a $(2^\kappa)^+$-saturated ideal on $\kappa$ was provable
in ZFC.  He focused on $\lambda$ with $\kappa\le\lambda\le 2^\kappa$.
Kunen pointed out that  arguments of Ulam \cite{Ulam1930} 
showed that if $\kappa$ is a successor cardinal, then there can be no 
$\lambda$-saturated ideal on $\kappa$ with $\omega<\lambda\le\kappa$.
He remarked (see \cite[p. 72]{Kunen1978})
that using the techniques of his 1970 paper \cite{Kunen1970},
from an $\omega_2$-saturated ideal on $\omega_1$ one gets consistency
of ``inner models with several measurable cardinals.''

\medskip

In 1974-1975, Alan Taylor\footnote{Based on email of May 8, 2016 from
  Alan Taylor.} and Stanley Wagon were both in Berkeley for nine months
and the Baumgartners spent the winter and spring quarters there.  Taylor
became interested in Wagon's work on saturation of ideals. 
At the end of 1975, Baumgartner, Taylor and Wagon \cite{JEB1977TW} submitted
their paper \emph{On splitting stationary subsets of large cardinals} which 
appeared in 1977.  They
looked at saturation properties of ideals, especially nonstationary ideals.
In 1972 Kunen \cite{Kunen1978} had shown the consistency of an 
$\omega_2$-saturated ideal on $\omega_1$ relative to the existence of
a huge cardinal, but only partial results were
available on when or if the non-stationary ideal on $\kappa$ could be 
$\kappa^+$-saturated.

Baumgartner, Taylor and Wagon \cite{JEB1977TW}
showed that given a normal ideal $\mathcal{I}$
on $\kappa$, $\mathcal{I}$ is $\kappa^+$-saturated if and only if the ideals
$\mathcal{I}\vert A$ generated by $I$ and $\kappa\setminus A$ for 
$A\in\mathcal{P}(\kappa)\setminus  \mathcal{I}$ are the only normal ideals 
that extend $\mathcal{I}$. Thus the non-stationary ideal NS${}_\kappa$ is 
$\kappa^+$-saturated  if and only if all normal ideals on $\kappa$ have the 
form NS${}_\kappa\vert A$ for some $A\subseteq \kappa$.  
As a corollary, they showed 
that if $\mathcal{I}$ is a normal $\kappa^+$-saturated ideal on $\kappa$ and $\mathcal{J}$ is a normal extension of $\mathcal{I}$, then $\mathcal{J}$ is also 
$\kappa^+$-saturated.  
It follows that if NS${}_\kappa$ is $\kappa^+$-saturated, then
every normal non-trivial ideal on $\kappa$ is  $\kappa^+$-saturated.
As a corollary, they showed that if $\kappa$ is 
greatly Mahlo, then the nonstationary ideal on $\kappa$, NS${}_\kappa$ is not 
$\kappa^+$-saturated.
Gitik and Shelah 
\cite{Gitik1997Shelah} proved that $\omega_1$ is the only uncountable 
cardinal for which NS${}_\kappa$ can be $\kappa^+$-saturated. 

Foreman \cite{Foreman2015} highlighted the power of non-stationary ideals
and their restrictions to selected stationary sets when he showed that
the consistency of ZFC together with his strengthening of the
classical Chang conjectures to the principle of \emph{Strong Chang
  Reflection}\footnote{We omit the definition of this principle but
  note that it includes second order reflection requirements and that
  Foreman has shown it is consistent from a $2$-huge cardinal.} for $(\omega_{n+3},\omega_n)$
implies the consistency of ZFC together with the existence of a huge
cardinal in a model of the form $L[A^\ast, \breve{I}]$ where $\breve{I}$ is
the dual of the appropriate nonstationary ideal.  Foreman used the
proposition below to show that the set $A^\ast$ was absolutely
definable.

\smallskip

\noindent{\bf Proposition} (Baumgartner): Let $M$, $N\prec
H(\theta)$.  
\newline
If 
$\sup(M\cap\omega_{n+2}) =  \sup(N\cap\omega_{n+2})\in\cof(>\omega)$,  
$N\cap\omega_{n+1}=M\cap\omega_{n+1}$ and
$\sup(N\cap\omega_{n+1})\in\cof(>\omega)$,
then $M\cap\omega_{n+2}=N\cap\omega_{n+2}$.

\medskip

In 1982, Baumgartner and Taylor \cite{JEB1982Tsatideals1},
\cite{JEB1982Tsatideals2} published a 
two part paper on saturation properties where they pioneered the study
of conditions on a forcing which preserved the saturation property of
ideals in the extension.
In the first part, given a cardinal $\lambda$, Baumgartner and Taylor 
concentrated on questions about which properties of an ideal $\mathcal{I}$ and 
a partial order $\mathbb{P}$ guarantee the $\lambda$-saturatedness of the ideal $\overline{\mathcal{I}}$ 
generated by $\mathcal{I}$ in the generic extension by $\mathbb{P}$.
They focus on instances that do not call for the use of large cardinals.  
For example they prove that 
if the forcing has
the $\sigma$-finite chain condition,\footnote{A partial order $P$
  satisfies the \emph{$\sigma$-finite chain condition} if there is a
  function $f:P\to\omega$ such that for all $n<\omega$, every pairwise
  incompatible subset of $f^{-1}(\{n\})$ is finite.}  
then one can conclude that in $M[G]$
every ideal on $\omega_1$ is $\omega_2$-generated, and hence, by a result
earlier in the paper, is $\omega_3$-saturated. 
They asked whether under ccc forcing, 
the converse that all $\omega_3$-saturated ideals are $\omega_2$-generated.
Baumgartner and Taylor used 
a ccc forcing GH which is a variant of one by Galvin and Hajnal and 
showed that in the extension, there is an ideal on $\omega_1$ that is not
$\omega_3$-saturated.  
In Corollary 3.5 they show that it is relatively
consistent with ZFC that $2^\omega$ is large and the nonstationary ideal 
$\mathcal{I}$
on $[\omega_2]^\omega$ is not $2^\omega$-saturated, but there is a stationary set
$S\subseteq [\omega_2]^\omega$ such that $\mathcal{I}\vert S$ is 
$\omega_4$-saturated.  

In part 2 
Baumgartner and Taylor 
\cite{JEB1982Tsatideals1} continue the study of preservation under forcing, 
especially ccc forcing, of saturation properties of countably complete ideals 
such as $\omega_2$-saturation or precipitousness.\footnote{Precipitous
  ideals were introduced by Jech and Prikry \cite{JechPrikry1976}.  If
    $I$ is an ideal on $\mathcal{P}(\kappa)$, then
    $\mathcal{P}(\kappa)/I$ is a notion of forcing which adds an
    ultrafilter $G$ extending the filter dual to $I$, and the ideal $I$ is
  said to be \emph{precipitous} if $\kappa\Vdash_{\mathcal{P}(I)/I} 
V^\kappa/G$ is wellfounded.} 
They formulate equivalences of the 
$\omega_2$-saturation of a countably complete ideal $\mathcal{I}$ on 
$\omega_1$ being preserved under ccc forcing in terms of a generalized 
version of Chang's conjecture and 
a weakening of Kurepa's Hypothesis.
They showed that given  an $\omega_2$-saturated ideal $\mathcal{I}$,
in any forcing extension by a $\sigma$-finite chain 
condition forcing, the ideal $\overline{\mathcal{I}}$ induced by
$\mathcal{I}$ is $\omega_2$-saturated on $\omega_1$.
They also showed that after forcing with the partial order for adding a closed
unbounded subset of $\omega_1$ with finite conditions (see \cite[page 926]
{JEB1984settop}), 
there are no $\omega_2$-saturated countably complete ideals on
$\omega_1$ in the extension.
They revisited the variant of the Galvin-Hajnal partial ordering GH used to
provide a consistent counterexample to all $\omega_2$-generated
countably complete ideals on $\omega_1$ being $\omega_3$-saturated,
and showed that the $\omega_2$-saturation of any ideal on $\omega_1$ 
is preserved when forcing with GH. 

Baumgartner and Taylor called an ideal $\mathcal{I}$  
\emph{presaturated} if it is both precipitous and $\omega_2$-preserving
i.e. $\Vdash_{\mathcal{P}(\omega_1)/\mathcal{I}}$ ``$\check{\omega}_2$ is a cardinal'').
After developing basic properties of presaturated ideals,  they prove that
(a) any $\omega_2$-preserving ideal on $\omega_1$ is a weak $p$-point; and 
(b) if there is a presaturated ideal on $\omega_1$, then there is a normal
presaturated ideal on $\omega_1$. 

A countably complete ideal $\mathcal{I}$ on $\omega_1$ is \emph{strong}
if and only if it is precipitous and
$\Vdash_{\mathcal{P}(\omega_1)/\mathcal{I}}{}j(\omega^V_1)=\omega^V_2$.
Baumgartner and Taylor observed that the argument given by 
Kunen \cite{Kunen1978} that the consistency of existence of a 
non-trivial countably complete $\omega_2$-saturated ideal 
on $\omega_1$ implies consistency of existence of several 
measurable cardinals, he only used the fact that the ideal was strong.  
Every $\omega_2$-saturated ideal is presaturated, and every
presaturated ideal is strong.

In section 5, Baumgartner and Taylor used a technical condition on a 
forcing $\mathbb{P}$, 
being \emph{$\mathcal{I}$-regular},\footnote{All 
ccc forcings and the forcing to add a closed unbounded subset of $\omega_1$
with finite conditions are $\mathcal{I}S$-regular.} to prove that if an ideal 
$\mathcal{I}$ in the ground 
model has one of the properties (a) precipitous, (b) strong, (c) presaturated, 
(d) $\omega_2$-saturated, then, in the extension by $\mathcal{P}$,
there is a set $A$ which is positive for the ideal  $\overline{\mathcal{I}}$ 
generated from $\mathcal{I}$ such that $\overline{\mathcal{I}}\vert A$ is 
has the corresponding property in the extension 
modulo the following constraints:
(a) precipitous (no additional constraint), 
(b) strong if $P$ does not collapse $\omega_2$, 
(c) presaturated if in the extension $J(\mathbb{\mathcal{P}})$ does not collapse
 $\omega_2^V$, 
(d) $\omega_2$-saturated if in the extension $j(\mathcal{P})$ is a ccc 
forcing). Moreover, if $\mathbb{P}$ is a ccc forcing then $A=\omega_1$ for
(a), (b), and (c).

In Theorem 5.10 they prove that the consistency of ZFC + ``there is a presaturated 
ideal on $\omega_1$'' implies the consistency of ZFC + ``there is a presaturated 
ideal on $\omega_1$ but no $\omega_2$-saturated ideals on $\omega_1$'' and
prove that consistency of ZFC + ``there is a precipitous 
ideal on $\omega_1$'' implies the consistency of ZFC + ``there is a precipitous
ideal on $\omega_1$ but no strong ideals on $\omega_1$'' (so also no
presaturated ideals on $\omega_1$).
%

We now review  a few of the  questions from the final section of the paper.
In Question 6.1, Baumgartner and Taylor \cite{JEB1982Tsatideals2} 
asked if the $\omega_2$-saturation of a countably complete ideal on 
$\omega_1$ is preserved under ccc forcing and in Question 6.2, they
asked if $\omega_1$ can carry a countable complete $\omega_2$-saturated ideal which
satisfies Chang's conjecture.  They had
proved for an $\omega_2$-saturated ideal $\mathcal{I}$, that
$\mathcal{I}$ satisfies Chang's conjecture if and only if for every ccc
partial ordering $\mathbb{P}$, $\Vdash_{\mathbb{P}}$ "$\mathcal{I}$
generates an  $\omega_2$-saturated ideal''. 
Foreman, Magidor and Shelah \cite[p. 24, Corollary 17]{Foreman1988MagShelah}
proved that if Martin's Maximum (MM) holds, then the NS${}_{\omega_1}$ is 
$\omega_2$-saturated and there is no ccc forcing which destroys its
saturation, so  NS${}_{\omega_1}$ satisfies Chang's conjecture, giving
a model in which the answers to Questions 6.1 and 6.2 are yes.
Donder and Levinski (unpublished) gave a model in which the answer is
no to Question 6.1.  
Boban Velickovic \cite{Velickovic1992forcestat} gave another negative answer
to Question 6.1 in a forcing extension of a model of MM.

In Question 6.5, Baumgartner and Taylor asked if 
every $\omega_2$-preserving countably complete
ideal on $\omega_1$ is precipitous. John Krueger 
\cite[page 844, Corollary 6]{Krueger2003} gave a positive answer for normal 
ideals on $\omega_1$ under the cardinality constraint $2^{\omega_1}\le \omega_3$
when he proved that if $\kappa$ is regular and $2^\kappa\le\kappa^{++}$, 
then the properties of $\kappa^+$-preserving and 
presaturated\footnote{What Baumgartner and Taylor called presaturated, 
Krueger called weakly presaturated, under the hypothesis of this
theorem, Kruegar proved the Baumgartner-Taylor version and his version were equivalent.} are equivalent for normal ideals.

In Question 6.10, Baumgartner and Taylor asked if the consistency of
the existence of a 
strong ideal on $\omega_1$ is equivalent to the consistency of the 
existence of a normal strong ideal. 
In 2010, Gitik \cite[page 196, Proposition 3.1]{Gitik2010normal} gave a 
positive answer. 

In Question 6.11, Baumgartner and Taylor asked if the consistency of the 
existence of a strong ideal on $\omega_1$ is equivalent to the consistency 
of the existence of an $\omega_2$-saturated ideal on $\omega_1$.
Consider the following four statements:
\begin{enumerate}
\item[(a)] ZFC + ``there exists a Woodin cardinal.''
\item[(b)] ZFC + ``there exists an $\omega_2$-saturated ideal on $\omega_1$.'' 
\item[(c)] ZFC + ``there exists a presaturated ideal on $\omega_1$.'' 
\item[(d)] ZFC + ``there exists a strong ideal on $\omega_1$.'' 
\end{enumerate}
That the consistency of (a) 
implies the consistency of (b) was shown
in a series of papers each building on the previous which dropped the large
cardinal  needed from a supercompact in 1983\footnote{See
  \cite{Foreman1988MagShelah}  for the proof. The timing is from a private communication from Foreman,
  November 18, 2016.}, to a Shelah
cardinal in 1984\footnote{See \cite{ShelahWoodin1990} for the proof
  and \cite{BagariaonShelahWoodin} for the timing of the major breakthrough.}
to a Woodin cardinal after its invention in 1984\footnote{See \cite{Steel2007} for the timing.} and no later than
1985\footnote{See  \cite{Shelah1987NSsat} for an announcement of the
  result and its timing, and 
  \cite{ShelahProperImproper} for a proof}. 
Baumgartner and Taylor \cite{JEB1982Tsatideals2} observed that an
$\omega_2$-saturated ideal is presaturated and strong,  so the
consistency of each statement implies the consistency of the next on
the list. John Steel and  Jensen 
\cite{JensenSteel2013} building on Steel \cite{Steel1996CoreIterate} 
proved the consistency of (c) implies the consistency of (a), 
and Benjamin Claverie and Ralf Schindler
\cite[Section 6]{ClaverieSchindler2012} proved the consistency of (d)
implies the consistency of (a). Hence all four statements are
equiconsistent and Question 6.11 is answered positively.

\subsection{Iterated forcing and Axiom A}
In 1976 Richard Laver \cite{Laver1976BorelConj} introduced the modern form of 
iterated countable support forcing in his celebrated paper on the consistency 
of the Borel Conjecture.  Other early countable support interations include
a term forcing of Mitchell \cite{Mitchell1970}, \cite{MitchellAtree1972}
and the forcing in Jensen's consistency proof of Suslin's Hypothesis with the 
Continuum Hypothesis \cite{DevlinJohns1974} which appeared in 1974.

In 1979 Baumgartner and Laver \cite{JEB1979Laver}
developed a countable support 
iterated Sacks forcing, and used it to prove the consistency 
of ZFC + $2^{\aleph_0}=\aleph_2$ with every selective ultrafilter\footnote{Selective ultrafilters are also known as Ramsey ultrafilters and as 
Rudin-Keisler minimal ultrafilters.} being $\aleph_1$-generated, solving
Erd\H{o}-Hajnal Problem 26 in  \cite{EH71.prob}.
They also used their forcing to give a new proof 
of the result of Mitchell \cite{Mitchell1972Atree}  that it is consistent that there are no
$\omega_2$-Aronszajn trees.

\medskip

Baumgartner's invention of Axiom~A forcing was a critical point in the 
development of generalizations of Martin's axiom.  A partial order $(\mathbb{P},\le)$ satisfies \emph{Axiom A} if and only if there is a sequence
$\langle \le_n: n\in\omega\rangle$ of partial orderings of $\mathbb{P}$ such that 
$p\le_0 q$ implies $p\le q$, for every $n$, $p\le_{n+1}q$ implies $p\le_n q$,
and the following conditions hold: 
\begin{enumerate}
\item if  $\langle p_n\in\mathbb{P}: n<\omega\rangle$ is a sequence
such that $p_0\ge_0 p_1\ge_1 \dots \ge_{n-1} p_n\ge_n \dots$ , then there is a 
$q\in\mathbb{P}$ such that $q\le_n p_n$ for all $n$;  
\item for every $p\in\mathbb{P}$, for every $n$ and for every ordinal name
$\dot{\alpha}$, there exist a $q\le_n p$ and a countable set $B$ such that
$q\Vdash\dot{\alpha}\in B$.
\end{enumerate}
The class of Axiom A forcings includes countable chain condition
forcings, countably closed forcings, Sacks (perfect set) forcing,
Prikry forcing, and Mathias forcing.  It is both a generalization of
ccc and $\sigma$-closed forcing.
Baumgartner proved the consistency of a forcing axiom generalizing Martin's axiom to cover all Axiom~A forcings, and 
in the summer of 1978 he included the proof in a series of lectures on iterated forcing  
in the three week long Summer School in Set Theory
in Cambridge, England organized by Harrington, Magidor and Mathias.
Baumgartner's expository paper \cite{JEB1983iterated}  growing out 
of these lectures was aimed at individuals with basic knowledge of forcing, and has been a popular introduction to the subject for graduate students 
for many 
years.\footnote{James Cummings \cite[page 7]{Cummings2010hand})  
describes Section 7 of his chapter \emph{Iterated Forcing and Elementary 
Embeddings} for the \emph{Handbook of Set Theory} ``as essentially following the 
approach of Baumgartner's survey.''}
Baumgartner indicated that the approach he took to iterated forcing
was ``strongly influenced by Laver's paper [11]''
\cite[page 2]{JEB1983iterated}\footnote{The paper of Laver [11] is 
his Borel Conjecture paper \cite{Laver1976BorelConj}} and thanked
Laver and Shelah for conversations and correspondence.

In 2005 Tetsuya Ishiu \cite{Ishiu2005AxiomA}  proved that a poset 
is forcing equivalent to a poset satisfying Axiom A if and only if 
it is $\alpha$-proper for every $\alpha<\omega_1$. 
A notion of forcing is \emph{proper} if for all regular 
uncountable cardinals $\lambda$, the forcing preserves stationary 
subsets of $[\lambda]^\omega$.  
Properness was developed by Shelah starting in 1978 and first
appeared in print in 1980 (see \cite{Shelah1980LinO}). 
Being $\alpha$-proper is a natural strengthening by Shelah of being
proper (see Shelah's book \cite{Shelah1982ProperForcing}).

\subsection{Proper forcing and the Proper Forcing Axiom}
Axiom~A forcing was in important influence in the development of 
the Proper Forcing Axiom (PFA), which further  extended the class of
applicable forcings to include proper forcing. 
Early in 1979 (see \cite[926]{JEB1984settop}), Baumgartner formulated
the Proper Forcing Axiom, which can
be briefly described as the extension of Martin's Axiom to include
proper forcing as stated below:
\begin{quote}
{\bf Proper Forcing Axiom (PFA): } 
If $\mathbb{P}$ is a proper forcing and $\mathcal{D}$ is a collection
of at most $\aleph_1$ dense sets, then there is a filter
$G\subseteq\mathbb{P}$ which meets every element of $\mathcal{D}$.  
\end{quote} 

Baumgartner used a Laver Diamond to prove that if ZFC together
with the existence of a supercompact cardinal.
is consistent, then so is ZFC + $2^{\aleph_0}=\aleph_2$ + PFA.
A   cardinal $\kappa$ is \emph{$\lambda$-supercompact} if there is an
  elementary embedding $j:V\to M$ so that the critical point of $j$ is
  $\kappa$, $j(\kappa)>\lambda$, and $M$ contains all its $\lambda$ sequences.
  A cardinal $\kappa$ is \emph{supercompact} if it is
  $\lambda$-supercompact for all $\lambda\ge\kappa$. Alternatively,
  $\kappa$ is \emph{supercompact} if for all $A$ of power at least
  $\kappa$, there is a normal measure on $[A]^{<\lambda}$.
A \emph{Laver diamond} \cite{Laver1978diamond} for a supercompact
cardinal $\kappa$ is a function $f:\kappa\to V_\kappa$ such that for
every $x\in V_\kappa$ and every $\lambda\ge |TC(x)$, there is a supercompact
ultrafilter $U_\lambda$ on $[\lambda]^{<\kappa}$ such that $(j_\lambda f)(\kappa)=x$.
The Laver Diamond was used to organize the critical iteration.

\medskip

In 1983, Baumgartner's expository paper on iterated forcing, which
was based his lectures at the 1978 Cambridge Summer School, appeared in 
the proceedings of that conference. Devlin published his 
\emph{Yorkshireman's Guide to Proper Forcing} 
\cite{Devlin1983Yguide} in the same proceedings in which he gave a proof of 
the Proper Forcing Axiom.  Devlin, who was not at the Cambridge
summer school,  was encouraged to write the article 
independently by Baumgartner, Rudi G\"obel, and Todorcevic. 
He based his article, which started out as personal notes,  on the following
materials:
\begin{itemize}
\item 
Notes written by Shelah in Berkeley in 1978 when he was giving lectures on proper forcing on material that eventually appeared in Chapters III, IV, V of the
first edition of \emph{Proper Forcing}; and
\item Notes written by Juris Steprans based on the ten lectures given by 
Baumgartner at the SETTOP Meeting in July and August, 1980 in Toronto.
\end{itemize}

\medskip

In 1984, the \emph{Handbook of Set-Theoretic Topology} was published and became
an important reference for set theorists and set-theoretic topologists.
Baumgartner \cite{JEB1984settop} wrote an extensive article, 
\emph{Applications of the Proper Forcing Axiom} starting from the definitions but
requiring knowledge of forcing.  He gave an example of a forcing that
was just barely proper, namely the forcing $P$ 
(see \cite[page 926]{JEB1984settop}) to add a club to $\omega_1$ with finite 
conditions. Conditions in $P$ are finite functions from $\omega_1$ into 
$\omega_1$ approximating an enumeration of a closed unbounded set, and
conditions are ordered by reverse inclusion.  This construction was
generalized by Todorcevic in \cite{Todor1984onPFA} in which he developed his 
seminal method for building proper partial orders using models as
side conditions.
Building on Baumgartner's elegant approach to adding a club to $\omega_1$,
Friedman \cite{Friedman2006}, Mitchell \cite{Mitchell2009},  and 
Neeman \cite{Neeman2014} all developed forcings with finite conditions to
add a club to $\omega_2$.  Inspired by the forcings of Friedman and Mitchell,
John Krueger defined \emph{adequate sets} and 
\emph{$S$-adequate sets}\footnote{We omit these definitions for brevity.} and 
developed a type of forcing for adding interesting combinatorial objects
with finite conditions using $S$-adequate sets of models as side conditions. 
In \cite{Krueger2015clubII} he used the approach to
add a closed, unbounded set to a given fat stationary set.

Baumgartner gave proofs from PFA 
of a number of statements known to be consistent some of which are listed below.
\begin{enumerate}
\item Theorem 6.9: PFA implies all $\aleph_1$-dense sets of reals are 
isomorphic.
\item Theorem 7.2: PFA implies there are no $\aleph_2$-Aronszajn trees. 
\item Theorem 7.10: PFA implies every tree of height $\omega_1$ and
  cardinality $\aleph_1$ is essentially special,
and therefore weak Kurepa's Hypothesis (wKH) is false.
\item Theorem 7.12: PFA implies that $\square_{\omega_1}$ is false.
\end{enumerate}

In Theorem 7.13 Baumgartner proved that PFA implies 
$\lozenge(E)$ for every stationary subset of
$\{\alpha<\omega_2: \cf\ \alpha=\omega_1\}$. 

Theorem 6.9 above is useful for the Basis Problem for uncountable linear orders.
As noted earlier, Baumgartner \cite{JEB1973reals}) proved
all $\aleph_1$-dense sets of reals are isomorphic relative to ZFC + 
$2^{\aleph_0}=\aleph_2$.  In 2006, Justin Moore \cite{MooreJT2006Basis}
used PFA to show there is a two element basis for the collection 
of Aronszajn lines, namely a Countryman type and its reverse answering
question 5.1(i) of \cite{JEB1976newotp}
Since the forcing Baumgartner  \cite{JEB1973reals}, \cite[Theorem 6.9]{JEB1984settop}
used in the consistency of a one element basis for the class of real types 
is proper, it can be combined with the consistency of a two element
basis for the class Aronszajn lines relative to PFA  to get a 
positive answer relative to PFA to Question 5(ii) of Baumgartner \cite{JEB1976newotp}, and
with the addition of $\omega_1$ and $\omega_1{}^\ast$, we obtain 
the consistency that the class of uncountable orderings has a five element 
basis relative to PFA. 

With regard to Theorem 7.2, Silver (see \cite{Mitchell1972Atree}) using a model of Mitchell
proved that the non-existence of an $\aleph_2$-Aronszajn tree is equiconsistent with the existence of a weakly compact cardinal.

With regard to Theorem 7.10, Mitchell \cite{Mitchell1972Atree} proved the 
failure of weak Kurepa's hypothesis is equiconsistent 
with the existence of an inaccessible cardinal over ZFC. 
Baumgartner \cite{JEB1983iterated} and independently and earlier Todorcevic 
\cite{Todor1981wKH} proved ZFC + MA + $\neg$wKH is consistent relative to 
the existence of an inaccessible cardinal, and Todorcevic gave consequences 
in his paper. 

Theorem 7.12 was improved by Todorcevic \cite{Todor1984onPFA}
who proved that PFA implies $\square_\kappa$ fails for all uncountable cardinals.
Note that the failure of $\square_\kappa$, for a regular uncountable cardinal $\kappa$,
is equiconsistent with the existence of a Mahlo cardinal.  In the early 1970s,
Solovay\footnote{For timing of the result, see 
Math Review MR2833150 (2012g:03134) by Kanamori; for the attribution see 
\cite{KruegerSchimmerling2011}, \cite[page 547]{Jech2003}.} 
proved that if $\lambda >\kappa$ is a Mahlo cardinal, then in
an extension by the L\'evy collapse Coll$(\kappa,<\lambda)$, 
$\lambda=\kappa^+$ and $\neg\square_\kappa$. 
Jensen \cite{Jensen1972fine} proved that if $\kappa^+$ is not Mahlo in $L$, 
then $\square_\kappa$ holds.

\smallskip

Baumgartner 
\cite[Section 8]{JEB1984settop} introduced a strengthening of PFA which he called PFA${}^+$ and under its assumption
proved two theorems on stationary-set  reflection and pointed
out that that such results imply the consistency of many measurable cardinals
(see \cite{Kanamori1978Magidor}).  

In 2009, Jensen, Schimmerling, Schindler and Steel \cite{StackingMice2009}
used the results of Todorcevic \cite{Todor1984onPFA}, that PFA implies $2^{\aleph_0}=\aleph_2$ and
$\square_\kappa$ fails for all uncountable cardinals, together with core model 
theory to show that PFA implies there 
is an inner model with a proper class of strong cardinals and a proper 
class of Woodin cardinals, and indiscernibles for such a model.

In 2011, Matteo Viale and Christoph Wei{\ss} 
\cite{Viale2011Weiss}
proved that if one can force PFA with a proper forcing that collapses a large 
cardinal $\kappa$ to $\omega_2$  and satisfies the $\kappa$-covering and
$\kappa$-approximation properties, then $\kappa$ is
supercompact. These papers 
of Jensen, Schimmerling, Schindler and Steel and of
Viale and Weiss suggest that Baumgartner's use of a supercompact cardinal in 
obtaining the consistency of PFA is likely necessary.

\subsection{Chromatic number of graphs}
Let us call a coloring of the vertices of a graph \emph{good}
if no pair joined by an edge have the same color.
an edge in the graph have the same color.
The chromatic number of a graph is the smallest number of colors 
for which there is a good coloring.
In 1984 Baumgartner \cite{JEB1984genericGraph} 
proved that If ZFC is consistent,
then so is ZFC $+$ GCH + ``there is a graph of cardinality $\aleph_2$ and 
chromatic number $\aleph_2$ such that every subgraph of cardinality 
$<\aleph_2$ has chromatic number $\le \aleph_0$'' providing 
a consistent negative answer to a question 
from 1961 of Erd\H{o}s and Hajnal \cite[page 118]{EP1961Hajnal} (this 
quote has been mildly rephrased with modern notation):
\begin{quote}
 Let there be given a graph $G$ of power $\aleph_2$. Suppose that
 every subgraph $G_1$ of $G$ of cardinality at most $\aleph_1$ has  
 chromatic number not greater than $\aleph_0$.  Is it then true 
that the chromatic number of $G$ is not greater than $\aleph_0$?
\end{quote}
The question was reiterated in print by Erd\H{o}s and Hajnal 
in 1966 in  \cite[pages 92-93]{EP1966Hajnal} in a paper  dedicated to the
to the 60th birthdays of well-known Hungarian mathematians R\'ozsa
P\'{e}ter and L\'aszl\'o Kalm\'ar both born in 1905.
In their 1968 paper emerging from a 1966 conference, 
Erd\H{o}s and Hajnal gave a negative answer under CH when they proved
there is a graph on $(2^{\aleph_0})^+$ many vertices whose chromatic
number is at least $\aleph_{1}$ and all of whose subgraphs of smaller 
cardinality have chromatic number at most $\aleph_0$. The statement of
the theorem was immediately followed by two questions:  (A) Does there  
exist a graph of power $\omega_{\omega+1}$ and uncountable chromatic
number all of whose smaller subgraphs have countable chromatic number?
(B) Does there exist a graph of power and chromatic number $\omega_2$
all of whose smaller subgraphs have countable chromatic number.
These appear as Problem 41(A) and 41(B) in the 
problem paper \cite{EH71.prob} growing out of
the presentation at the 1967 UCLA summer school and were reiterated  
in 1975 in \cite[page 415]{EP1975combo.prob} and in 1973
and Galvin \cite{Galvin1973chrom} reframed one question by asking if every
graph of chromatic number $\aleph_2$ has a subgraph of chromatic number
$\aleph_1$.   Thus Baumgartner gave a positive answer to problem 41(B)
and a negative answer to Galvin's problem. In 1988, Komjath%
\footnote{In 2002, Komjath \cite{Komjath2002}
made a systematic study of the set of
of chromatic numbers realized by subgraphs of a given graph.}
\cite{Komjath1988} provided a different positive answer to problem 41(B)
since in his forcing extension $2^{\aleph_0}=\aleph_3$.
Also in 1988, Foreman and Laver \cite{ForemanLaver1988} proved the 
relative consistency of the opposite conclusions, assuming the 
existence of a huge cardinal to construct  a forcing extension in 
which ZFC + GCH hold and every graph of power $\aleph_2$ and chromatic
number $\aleph_2$ has a subgraph of size and chromatic number $\aleph_1$. 
In 1997 Todorcevic \cite{Todor1997}  constructed in ZFC a graph
of power $2^{\aleph_2}$ of uncountable chromatic number with no
subgraph of power and chromatic number $\aleph_1$.
Recent work on the construction of graphs of large power 
and  uncountable chromatic number all of whose smaller subgraphs
have countable chromatic number includes  
the use of $\square_\lambda$ + $2^{\lambda}=\lambda^+$ for an uncountable
cardinal $\lambda$ by Assaf Rinot \cite{Rinot2015} to get graphs of size
$\lambda^+$ of chromatic number of any desired value $\kappa\le
\lambda$. 

\subsection{A thin very-tall superatomic Boolean algebra}
In 1987  Baumgartner and Shelah published their proof of
the consistent existence of  a thin very-tall superatomic Boolean
algebra, where a 
Boolean algebra is superatomic  if and only all of its homomorphic
images are atomic.  Their collaboration came about as follows.
Baumgartner circulated a preprint which included his proof of this result
by a two step forcing, a countably closed forcing and a ccc forcing. Baumgartner used 
a function $f_*$ with special properties in his construction of the second 
forcing to guarantee it was ccc. Fleissner found an error in the proof and it was later discovered that no function exists with the special properties Baumgartner had envisioned. Shelah came up with a different set of properties $\Delta$ 
for a function $f_*$ and a different countably closed forcing to make the whole 
construction work. The combined forcing of Baumgartner and Shelah \cite{JEB1987Shelah} 
has proven useful in other settings.

In 2001, Juan Carlos Mart{\'{\i}}nez \cite{Martinez2001} generalized the
result to show that Con(ZFC) implies Con(ZFC + ``for all
$\alpha<\omega_3$, there is a superatomic Boolean algebra of width $\omega$ and
height $\alpha$'').

In 2013, Boban Velickovic and Giorgio Venturi \cite{VelickovicVenturi}
used Neeman's method of forcing 
with generalized side conditions with two types of models and finite support
to give a new proof of the Baumgartner-Shelah result.

Recent research in the area (see \cite{Martinez2016}) has turned to proofs of existence of cardinal
sequences of locally compact scattered (Hausdorff) spaces (LCS), and
these results can be translated into results about superatomic Boolean algebras.

\subsection{Closed unbounded sets}
In 1991, Baumgartner \cite{JEB1991} published the paper on the structure
of closed unbounded subsets of $[\lambda]^{<\kappa}$, and the stationary
sets associated with them, with a focus on $\lambda=\kappa^{+n}$.
Given cardinals 
$\kappa<\lambda$, he introduced a family of sets 
$S(\kappa,\lambda;\kappa_0, \dots , \kappa_n)$ 
parameterized for some $n<\omega$ by a sequence of regular cardinals 
$\kappa_0,\kappa_1,\dots ,\kappa_n$ all smaller than $\kappa$, and
he proved that they are stationary sets.
His stated goal was to find closed unbounded sets $C$ so that 
each $x\in C\cap 
S(\kappa,\lambda;\kappa_0, \dots , \kappa_n)$  was determined as much as 
possible by the sequence 
$\left\langle 
\sup (x\cap \kappa_0), \ldots , \sup (x\cap \kappa_n)\right\rangle$.
He used these sets and their intersections to prove limitations on the
size of intersections of closed unbounded sets with these 
sets $S(\kappa,\lambda;\kappa_0, \dots , \kappa_n)$.  
For example, for $\kappa$ regular and $\lambda=\kappa^{+n}$, Baumgartner
proved that if all the $\kappa_i$'s are regular,  
$\kappa_i=0$ for some $i>0$, and $\lambda=\kappa^{+n}$, 
then the intersection of every
closed unbounded set $C$ with 
$S(\kappa,\lambda;\kappa_0, \dots , \kappa_n)$  has cardinality at least $\lambda^\omega$.
 
He also introduced two combinatorial principles which were useful in pinning
down the intersections of closed unbounded sets with stationary sets of the
form $S(\kappa,\lambda;\kappa_0, \dots , \kappa_n)$.
Then he used various types of
$\square$-sequences, the  weakening of Erd\H{o}s cardinals to remarkable 
cardinals, and notions of reverse-Easton-like forcings to prove consistency 
with and independence from ZFC of the combinatorial principles.

For a cardinal $\kappa$ he wrote $\square(\kappa)$%
\footnote{
Currently
  $\square(\kappa)$ denotes 
 a sequence $\langle   C_\alpha\mid\alpha<\kappa\rangle$
 such that $C_\alpha$ for limit $\alpha$ is a closed unbounded subset
  of $\alpha$; if $\beta$ is a limit point of $C_\alpha$, then
  $\beta\cap C\alpha=C_\beta$; and there is no \emph{threading},
  i.e. no closed unbounded set $C\subseteq\kappa$ such that
  $C\cap\alpha=C_\alpha$ for all limit points $C$.  This variant of Jensen's $\square_\lambda$
  principle is due to Todorcevic \cite{Todor1987partPairs}. See \cite[page 298]{Rinot2014} for context.%
}
to indicate that there is $\langle
C_\alpha\mid\alpha<\kappa,\mbox{ $\kappa$ singular limit}\rangle$ such
that for each singular limit $\alpha<\kappa$, $C_\alpha$ is closed and
unbounded subset of $\alpha$ of order type $<\alpha$; and if
$\beta<\alpha$ is a limit point of $C_\alpha$, then $C_\beta=\beta\cap
C_\alpha$ and indicated that $\square(\kappa^+)$ was equivalent to
$\square_\kappa$.  
In 1980 it was shown in \cite{BellerLitman1980} that if
$V=L$, then $\square(\kappa)$ holds for all $\kappa$.
In Theorem 6.6 Baumgartner proved that if $\kappa$ is $\gamma$-Erd\H{o}s for
some $\gamma<\kappa$, then there is a forcing extension in which
$\kappa$ remains $\gamma$-Erd\H{o}s and both $\square(\kappa)$ and
$\square(\{\alpha<\kappa\mid \alpha\mbox{ is regular}\})$ hold.

Foreman and Magidor \cite[page 66, Corollary 2.11(b)]{Foreman1995Magidor} 
generalized the function $S(\kappa,\lambda;\kappa_0,\kappa_1,\kappa_2,\kappa_3)$
by replacing $\lambda$ with the collection $H(\lambda)$
of sets hereditarily of cardinality $<\lambda$). 
While Baumgartner used the sets $S(\kappa,\lambda;\kappa_0, \dots , \kappa_n)$
to control sizes of clubs and their intersections with the stationary set,
Foreman and Magidor use their version to construct interesting stationary sets,
e.g. there is a non-reflecting stationary subset $C\subseteq 
S(\omega_2,H(\lambda); \omega_2,\omega_3; \omega_1,\omega)$.

\subsection{Revisiting partition relations}
A key result of Erd\H{o}s and Rado was the following theorem 
\cite[pages 467-8]{EP1956Rado}. 
\begin{quote}
{\bf Positive Stepping Up Lemma} (modern form): For all infinite cardinals
$\kappa$, all $\gamma$ with $2\le\gamma<\kappa$, all finite $r$, and
all cardinals $\langle \alpha_\nu:\nu<\gamma\rangle$, 
if $\kappa\rightarrow (\alpha_\nu)^r_{\nu<\gamma}$, 
then $(2^{<\kappa})^+ \rightarrow (\alpha_\nu+1)^{r+1}_{\nu<\gamma}$
\end{quote}
Let $\exp_n(\kappa)$ denote $n$-times iterated exponentiation,
that is, $\exp_0(\kappa)=\kappa$ and
$\exp_{n+1}(\kappa)=2^{\,\exp_n(\kappa)}$.  With this notation
and the cardinal arithmetic above, we state below the 
modern version of their Theorem 39, obtained
using the Positive Stepping Up Lemma repeatedly and starting
from the clear fact that $\exp_0(\kappa)=\kappa\rightarrow(\kappa)^1_\gamma$
for $\gamma<\cf(\kappa)$.  
\begin{quote}
\noindent{\bf Erd\H{o}s-Rado Theorem} 
(modern form): For every infinite cardinal $\kappa$,
every\index{Erd\H{o}s, P.}\index{Rado, R.} 
finite $r\ge 2$, and all $\gamma<\cf(\kappa)$,
$\left(\exp_r(2^{<\kappa})\right)^+ \rightarrow (\kappa+(r-1))^{r}_\gamma$.
\end{quote}

\smallskip

In 1993,
Baumgartner, Hajnal and Todorcevic \cite{JEB1993HT} 
published their extensions of 
the Erd\H{o}s-Rado Theorem
as follows.   

\begin{theorem}[Balanced Baumgartner-Hajnal-Todorcevic Theorem]
Suppose $\kappa$ is regular and uncountable.  For all $\ell<\omega$ and all 
ordinals $\xi$ with $2^\xi<\kappa$, 
$(2^{<\kappa})^+\rightarrow (\kappa+\xi)^2_\ell$.
\end{theorem}

\begin{theorem}[Unbalanced Baumgartner-Hajnal-Todorcevic Theorem]
Suppose $\kappa$ is regular and uncountable.  For all $\ell, n<\omega$, 
\[(2^{<\kappa})^+\rightarrow (\kappa^{\omega+2}+1,\kappa+n)^2.\]
\end{theorem} 

In 2001, Baumgartner and Hajnal \cite{JEB2001Hajnal} published
a polarized version of the Erd\H{o}s-Rado Theorem for pairs.
Baumgartner and Hajnal proved that
for every cardinal $\kappa$, 
\[\left(\begin{array}{c}(2^{<\kappa})^{++}\\ (2^{<\kappa})^+\end{array}\right)
\rightarrow 
\left(\begin{array}{c}\kappa\\ \kappa\end{array}\right)_{<\operatorname{cf}\kappa} .\]
If $\kappa$ is weakly compact, then this can be improved to
$\left(\begin{array}{c}\kappa^+\\ \kappa\end{array}\right) 
\rightarrow
\left(\begin{array}{c}\kappa\\ \kappa+1\end{array}\right)_{<\kappa}$ .

In 2003, Matthew Foreman and Hajnal \cite{Foreman2003Hajnal} used 
techniques from  \cite{JEB1993HT} and ideals  
coming from elementary submodels in their proof that if $\kappa^{<\kappa}$ 
and $\kappa$ carries a $\kappa$-dense ideal,
then $\kappa^+\rightarrow (\kappa^2+1,\alpha)^2$ for all $\alpha<\kappa^+$.
A proper, non-principal ideal $\mathcal{J}$ on $\kappa$ is \emph{$\lambda$-dense} if the 
Boolean algebra
$\mathcal{P}(\kappa)/\mathcal{J}$ has a dense set of size $\lambda$.
Equivalently, $\mathcal{J}$ is \emph{$\lambda$-dense} if
there is a family $\mathcal{D}$ of $\lambda$ many 
$\mathcal{J}$-positive sets so that every $\mathcal{J}$-positive set is 
$J$-almost included in some element of $\mathcal{D}$ and any two different elements of $\mathcal{D}$ have intersection in $\mathcal{J}$.  

Moreover, if $\kappa$ is a measurable cardinal, then there is a 
rather large ordinal $\Omega<\kappa^+$ such that for all $n<\omega$ and
$\alpha<\Omega$, $\kappa^+\rightarrow (\alpha)^2_n$.  They give multiple
definitions of $\Omega$ and use them to show that $\Omega$ is rather large. 
Specifically, they observe that it follows from the definitions  that
$L_\Omega$ is a model of ZFC; they show that the statement $``\alpha<\Omega$'' is
upwards absolute; and for $U$ a normal ultrafilter on $\kappa$,
they show that the least ordinal $\nu$ such that $L_\nu[U]\cap\kappa^{<\kappa}
=L[U]\cap\kappa^{<\kappa}$ is a lower bound.

In 2006, Albin Jones \cite{Jones2006} extended the weakly compact polarized 
partition relation of Baumgartner and Hajnal \cite{JEB2001Hajnal} to 
specify the order type of sets chosen to be homogeneous.

In 2014, Ari Brodsky gave a generalization of the
Balanced Baumgartner-Hajnal-Todorcevic Theorem to all partially ordered sets.
More specifically, Brodsky \cite{Brodsky2014} proved that if $P$
is a partial order such that $P\rightarrow
(2^{<\kappa})^1_{2^{<\kappa}}$ for some uncountable regular cardinal
$\kappa$, and if $\ell<\omega$ and $\xi$ is an ordinal such that
$2^{| \xi |}<\kappa$, then $P\rightarrow(\kappa+\xi)^2_\ell$.

\bigskip

\noindent{\bf Acknowledgements: }  
My appreciation to all the people with whom I have
  had conversations about Baumgartner in recent years including
  Yolanda Baumgartner, Joan Bagaria, Andres Caicedo, James Cummings, Natasha Dobrinen, Mirna
  D\v{z}jamonja, Monroe Eskew, Matt
  Foreman, Marcia Groszek, Andr\'as Hajnal, Akihiro Kanamori, Albin
  Jones, P\'eter Komj\'ath, 
Paul Larson,   Adrian Mathias, Bill Mitchell, Justin Moore, Assaf Rinot, Saharon Shelah, Frank
  Tall, Alan Taylor, and
  Stevo Todorcevic.  My apologies to anyone I have forgotten to mention.

\nocite{JEB1968dense.abs}
\nocite{JEB1969undefinable.abs} 
\nocite{JEB1969undefinable.abs}
\nocite{JEB1970decomp.abs} 
\nocite{JEB1971CantorThm.abs}
\nocite{JEB1974.ER}
\nocite{JEB1974topology.abs} 
\nocite{JEB1974onHindman}
\nocite{JEB1974topology.abs}
\nocite{JEB1974partition.abs}
\nocite{JEB1975vector}
\nocite{JEB1975.unctble}
\nocite{JEB1975HMsaturation}
\nocite{JEB1975GLMcKgame}
\nocite{JEB1975.unctble}
\nocite{JEB1975canon}
\nocite{JEB1978Galvin}
\nocite{JEB1978indy.abs}
\nocite{JEB1978TW}
\nocite{JEB1978Taylor}
\nocite{JEB1978indy.abs} 
\nocite{JEB1979survey}
\nocite{JEB1981K} 
\nocite{JEB.Boole1980} 
\nocite{JEB1982Weese}
\nocite{JEB1984genericGraph}
\nocite{JEB1984EPHiggs}
\nocite{JEB1984EPHiggs} 
\nocite{JEB1985SacksvsMA} 
\nocite{JEB1985Dordal}
\nocite{JEB1985base} 
\nocite{JEB1987Hajnal}
\nocite{JEB1988-9BalcerzakHejduk} 
\nocite{JEB1989partOrd}
\nocite{JEB1989polar}
\nocite{JEB1990Larson}
\nocite{JEB1990FZ} 
\nocite{JEB1990vD} 
\nocite{JEB1990} 
\nocite{JEB1991Spinas}
\nocite{JEB1993ST}
\nocite{JEB1994future}
\nocite{JEB1997onEP}
\nocite{JEB1995ufonN} 
\nocite{JEB2002Hajnal} 
\nocite{JEB2002Tall}

\bibliographystyle{plain}%
\bibliography{jebfull,near_jeb}

\begin{thebibliography}{100}

\bibitem{AbrahamRubinShelah1985}
Uri Abraham, Matatyahu Rubin, and Saharon Shelah.
\newblock On the consistency of some partition theorems for continuous
  colorings, and the structure of {$\aleph_1$}-dense real order types.
\newblock {\em Ann. Pure Appl. Logic}, 29(2):123--206, 1985.

\bibitem{AbrahamShelah1983}
Uri Abraham and Saharon Shelah.
\newblock Forcing closed unbounded sets.
\newblock {\em J. Symbolic Logic}, 48(3):643--657, 1983.

\bibitem{BagariaonShelahWoodin}
Joan Bagaria.
\newblock Book {R}eview: {L}arge cardinals imply that every reasonably
  definable set of reals is {L}ebesgue measurable by {S}aharon {S}helah and
  {H}ugh {W}oodin.
\newblock {\em Bulletin of Symbolic Logic}, 8(8):543--545, 2002.

\bibitem{JEB1988-9BalcerzakHejduk}
Marek Balcerzak, Jacek Hejduk, and James~E. Baumgartner.
\newblock On certain {$\sigma$}-ideals of sets.
\newblock {\em Real Anal. Exchange}, 14(2):447--453, 1988/89.

\bibitem{Banach1930}
Stefan Banach.
\newblock Uber additive massfunktionen in alstrakten mengen.
\newblock {\em Fund. Math.}, 15:97--101, 1930.

\bibitem{JEB1968dense.abs}
James~E. Baumgartner.
\newblock On the cardinality of dense subsets of linear orderings {I}.
\newblock {\em Notices Amer. Math. Soc.}, 15:935, 1968.
\newblock Abstract, preliminary report.

\bibitem{JEB1969undefinable.abs}
James~E. Baumgartner.
\newblock Undefinability of $n$-ary relations from unary functions.
\newblock {\em Notices Amer. Math. Soc.}, 17:842--843, 1969.
\newblock Abstract, preliminary report.

\bibitem{JEB1970decomp.abs}
James~E. Baumgartner.
\newblock Decompositions and embeddings of trees.
\newblock {\em Notices Amer. Math. Soc.}, 18:967, 1970.
\newblock Abstract, preliminary report.

\bibitem{JEB1970thesis}
James~E. Baumgartner.
\newblock {\em Results and independence proofs in combinatorial set theory}.
\newblock PhD thesis, University of California, Berkeley, 1970.

\bibitem{JEB1971CantorThm.abs}
James~E. Baumgartner.
\newblock A possible extension of {C}antor's theorem on the rationals.
\newblock {\em Notices Amer. Math. Soc.}, 18:428--429, 1971.
\newblock Abstract, preliminary report.

\bibitem{JEB1973reals}
James~E. Baumgartner.
\newblock All {$\aleph \sb{1}$}-dense sets of reals can be isomorphic.
\newblock {\em Fund. Math.}, 79(2):101--106, 1973.

\bibitem{JEB1974Hanf}
James~E. Baumgartner.
\newblock The {H}anf number for complete ${L}_{\omega_1,\omega}$-sentences
  (without {GCH}).
\newblock {\em J. Symbolic Logic}, 39:575--578, 1974.

\bibitem{JEB1974.ER}
James~E. Baumgartner.
\newblock Improvement of a partition theorem of {E}rd{\H{o}}s and {R}ado.
\newblock {\em J. Combin. Theory Ser. A}, 17:134--137, 1974.

\bibitem{JEB1974onHindman}
James~E. Baumgartner.
\newblock A short proof of {H}indman's theorem.
\newblock {\em J. Combinatorial Theory Ser. A}, 17:384--386, 1974.

\bibitem{JEB1974partition.abs}
James~E. Baumgartner.
\newblock Some results in the partition calculus.
\newblock {\em Notices Amer. Math. Soc.}, 21:A--29, 1974.
\newblock Abstract, preliminary report.

\bibitem{JEB1975canon}
James~E. Baumgartner.
\newblock Canonical partition relations.
\newblock {\em J. Symbolic Logic}, 40(4):541--554, 1975.

\bibitem{JEB1975.ineff}
James~E. Baumgartner.
\newblock Ineffability properties of cardinals. {I}.
\newblock In {\em Infinite and {F}inite {S}ets ({C}olloq., {K}eszthely, 1973;
  dedicated to {P}.~{E}rd\H os on his 60th birthday), {V}ol. {I}}, volume~10 of
  {\em Colloq. Math. Soc. J{\'a}nos Bolyai}, pages 109--130. North-Holland,
  Amsterdam, 1975.

\bibitem{JEB1975.unctble}
James~E. Baumgartner.
\newblock Partition relations for uncountable ordinals.
\newblock {\em Israel J. Math.}, 21(4):296--307, 1975.

\bibitem{JEB1975vector}
James~E. Baumgartner.
\newblock Partitioning vector spaces.
\newblock {\em J. Combinatorial Theory Ser. A}, 18:231--233, 1975.

\bibitem{JEB1974topology.abs}
James~E. Baumgartner.
\newblock Topological properties of {S}pecker types.
\newblock {\em Notices Amer. Math. Soc.}, 22:A--219, 1975.
\newblock Abstract, preliminary report.

\bibitem{JEB1976ad}
James~E. Baumgartner.
\newblock Almost-disjoint sets, the dense set problem and the partition
  calculus.
\newblock {\em Ann. Math. Logic}, 9(4):401--439, 1976.

\bibitem{JEB1976newotp}
James~E. Baumgartner.
\newblock A new class of order types.
\newblock {\em Ann. Math. Logic}, 9(3):187--222, 1976.

\bibitem{JEB1977ineff2}
James~E. Baumgartner.
\newblock Ineffability properties of cardinals. {II}.
\newblock In {\em Logic, foundations of mathematics and computability theory
  ({P}roc. {F}ifth {I}nternat. {C}ongr. {L}ogic, {M}ethodology and {P}hilos. of
  {S}ci., {U}niv. {W}estern {O}ntario, {L}ondon, {O}nt., 1975), {P}art {I}},
  pages 87--106. Univ. Western Ontario Ser. Philos. Sci., Vol. 9. Reidel,
  Dordrecht, 1977.

\bibitem{JEB1978indy.abs}
James~E. Baumgartner.
\newblock Independence results in set theory.
\newblock {\em Notices Amer. Math. Soc.}, 25:A248--249, 1978.

\bibitem{JEB1979survey}
James~E. Baumgartner.
\newblock Independence proofs and combinatorics.
\newblock In {\em Relations between combinatorics and other parts of
  mathematics ({P}roc. {S}ympos. {P}ure {M}ath., {O}hio {S}tate {U}niv.,
  {C}olumbus, {O}hio, 1978)}, Proc. Sympos. Pure Math., XXXIV, pages 35--46.
  Amer. Math. Soc., Providence, R.I., 1979.

\bibitem{JEB.Boole1980}
James~E. Baumgartner.
\newblock Chains and antichains in {${\cal P}(\omega )$}.
\newblock {\em J. Symbolic Logic}, 45(1):85--92, 1980.

\bibitem{JEB1982order}
James~E. Baumgartner.
\newblock Order types of real numbers and other uncountable orderings.
\newblock In {\em Ordered {S}ets ({B}anff, {A}lta., 1981)}, volume~83 of {\em
  NATO Adv. Study Inst. Ser. C: Math. Phys. Sci.}, pages 239--277. Reidel,
  Dordrecht, 1982.

\bibitem{JEB1983iterated}
James~E. Baumgartner.
\newblock Iterated forcing.
\newblock In {\em Surveys in set theory}, volume~87 of {\em London Math. Soc.
  Lecture Note Ser.}, pages 1--59. Cambridge Univ. Press, Cambridge, 1983.

\bibitem{JEB1984settop}
James~E. Baumgartner.
\newblock Applications of the {P}roper {F}orcing {A}xiom.
\newblock In {\em Handbook of {S}et-{T}heoretic {T}opology}, pages 913--959.
  North-Holland, Amsterdam, 1984.

\bibitem{JEB1984genericGraph}
James~E. Baumgartner.
\newblock Generic graph construction.
\newblock {\em J. Symbolic Logic}, 49(1):234--240, 1984.

\bibitem{JEB1985base}
James~E. Baumgartner.
\newblock Bases for {A}ronszajn trees.
\newblock {\em Tsukuba J. Math.}, 9(1):31--40, 1985.

\bibitem{JEB1985SacksvsMA}
James~E. Baumgartner.
\newblock Sacks forcing and the total failure of {M}artin's {A}xiom.
\newblock {\em Topology Appl.}, 19(3):211--225, 1985.

\bibitem{JEB1986onKunen}
James~E. Baumgartner.
\newblock Book {R}eview: {S}et theory. {A}n introduction to independence proofs
  by {K}enneth {K}unen.
\newblock {\em J. Symbolic Logic}, 51(2):462--464, 1986.

\bibitem{JEB1989polar}
James~E. Baumgartner.
\newblock Polarized partition relations and almost-disjoint functions.
\newblock In {\em Logic, {M}ethodology and {P}hilosophy of {S}cience, {VIII}
  ({M}oscow, 1987)}, pages 213--222. North-Holland, Amsterdam, 1989.

\bibitem{JEB1989partOrd}
James~E. Baumgartner.
\newblock Remarks on partition ordinals.
\newblock In {\em Set theory and its applications ({T}oronto, {ON}, 1987)},
  volume 1401 of {\em Lecture Notes in Math.}, pages 5--17. Springer, Berlin,
  1989.

\bibitem{JEB1990}
James~E. Baumgartner.
\newblock Is there a different proof of the {E}rd{\H{o}}s-{R}ado theorem?
\newblock In {\em A tribute to Paul Erd\H{o}s}, pages 27--37. Cambridge Univ.
  Press, Cambridge, 1990.

\bibitem{JEB1991}
James~E. Baumgartner.
\newblock On the size of closed unbounded sets.
\newblock {\em Ann. Pure Appl. Logic}, 54(3):195--227, 1991.

\bibitem{JEB1994future}
James~E. Baumgartner.
\newblock The future of modern set theory.
\newblock {\em Ann. Japan Assoc. Philos. Sci.}, 8(4):187--190, 1994.

\bibitem{JEB1995ufonN}
James~E. Baumgartner.
\newblock Ultrafilters on {$\omega$}.
\newblock {\em J. Symbolic Logic}, 60(2):624--639, 1995.

\bibitem{JEB1997onEP}
James~E. Baumgartner.
\newblock In {M}emoriam: {P}aul {E}rd{\H{o}}s 1913-1996.
\newblock {\em B. Symbolic Logic}, 3(1):70--71, 1997.

\bibitem{JEB2002Hajnal}
James~E. Baumgartner.
\newblock Hajnal's contributions to combinatorial set theory and the partition
  calculus.
\newblock In {\em Set {T}heory ({P}iscataway, {NJ}, 1999)}, volume~58 of {\em
  DIMACS Ser. Discrete Math. Theoret. Comput. Sci.}, pages 25--30. Amer. Math.
  Soc., Providence, RI, 2002.

\bibitem{JEB1985Dordal}
James~E. Baumgartner and Peter Dordal.
\newblock Adjoining dominating functions.
\newblock {\em J. Symbolic Logic}, 50(1):94--101, 1985.

\bibitem{JEB1984EPHiggs}
James~E. Baumgartner, Paul Erd{\H{o}}s, and Denis~A. Higgs.
\newblock Cross-cuts in the power set of an infinite set.
\newblock {\em Order}, 1(2):139--145, 1984.

\bibitem{JEB1990FZ}
James~E. Baumgartner, Ryszard Frankiewicz, and Pawel Zbierski.
\newblock Embedding of {B}oolean algebras in {$P(\omega)/{\rm fin}$}.
\newblock {\em Fund. Math.}, 136(3):187--192, 1990.

\bibitem{JEB1978Galvin}
James~E. Baumgartner and Fred Galvin.
\newblock Generalized {E}rd{\H o}s cardinals and {$0^{\#}$}.
\newblock {\em Ann. Math. Logic}, 15(3):289--313 (1979), 1978.

\bibitem{JEB1975GLMcKgame}
James~E. Baumgartner, Fred Galvin, Richard Laver, and Ralph McKenzie.
\newblock Game theoretic versions of partition relations.
\newblock In {\em Infinite and {F}inite {S}ets ({C}olloq., {K}eszthely, 1973;
  dedicated to {P}.~{E}rd\H os on his 60th birthday), {V}ol. {I}}, volume~10 of
  {\em Colloq. Math. Soc. J{\'a}nos Bolyai}, pages 131--135. North-Holland,
  Amsterdam, 1975.

\bibitem{JEB1973Hajnal}
James~E. Baumgartner and Andr{\'a}s Hajnal.
\newblock A proof (involving {M}artin's {A}xiom) of a partition relation.
  {Y}ellow {S}eries of the {U}niversity of {C}algary, {R}esearch {P}aper 122,
  {A}pril, 1971.
\newblock {\em Fund. Math.}, 78(3):193--203, 1973.

\bibitem{JEB1987Hajnal}
James~E. Baumgartner and Andr{\'a}s Hajnal.
\newblock A remark on partition relations for infinite ordinals with an
  application to finite combinatorics.
\newblock In {\em Logic and {C}ombinatorics ({A}rcata, {C}alif., 1985)},
  volume~65 of {\em Contemp. Math.}, pages 157--167. Amer. Math. Soc.,
  Providence, RI, 1987.

\bibitem{JEB2001Hajnal}
James~E. Baumgartner and Andr{\'a}s Hajnal.
\newblock Polarized partition relations.
\newblock {\em J. Symbolic Logic}, 66(2):811--821, 2001.

\bibitem{JEB1975HMsaturation}
James~E. Baumgartner, Andr{\'a}s Hajnal, and Attila M{\'a}t{\'e}.
\newblock Weak saturation properties of ideals.
\newblock In {\em Infinite and {F}inite {S}ets ({C}olloq., {K}eszthely, 1973;
  dedicated to {P}.~{E}rd\H os on his 60th birthday), {V}ol. {I}}, volume~10 of
  {\em Colloq. Math. Soc. J{\'a}nos Bolyai}, pages 137--158. North-Holland,
  Amsterdam, 1975.

\bibitem{JEB1993HT}
James~E. Baumgartner, Andr{\'a}s Hajnal, and Stevo Todorcevic.
\newblock Extensions of the {E}rd{\H{o}}s-{R}ado theorem.
\newblock In {\em Finite and {I}nfinite {C}ombinatorics in {S}ets and {L}ogic
  ({B}anff, {AB}, 1991)}, pages 1--17. Kluwer Acad. Publ., Dordrecht, 1993.

\bibitem{JEB1976HK}
James~E. Baumgartner, Leo Harrington, and Eugene Kleinberg.
\newblock Adding a closed unbounded set.
\newblock {\em J. Symbolic Logic}, 41:481--482, 1976.

\bibitem{JEB1981K}
James~E. Baumgartner and P{\'eter}. Komj{\'a}th.
\newblock Boolean algebras in which every chain and antichain is countable.
\newblock {\em Fund. Math.}, 111(2):125--133, 1981.

\bibitem{JEB1990Larson}
James~E. Baumgartner and Jean~A. Larson.
\newblock A diamond example of an ordinal graph with no infinite paths.
\newblock {\em Ann. Pure Appl. Logic}, 47(1):1--10, 1990.

\bibitem{JEB1979Laver}
James~E. Baumgartner and Richard Laver.
\newblock Iterated perfect-set forcing.
\newblock {\em Ann. Math. Logic}, 17(3):271--288, 1979.

\bibitem{BaxiomaticSetTheory1984}
James~E. Baumgartner, Donald~A. Martin, and Saharon Shelah, editors.
\newblock {\em Axiomatic set theory}, volume~31 of {\em Contemporary
  Mathematics}, Providence, RI, 1984. American Mathematical Society.

\bibitem{JEB1976Prikry}
James~E. Baumgartner and Karel Prikry.
\newblock On a theorem of {S}ilver.
\newblock {\em Discrete Math.}, 14(1):17--21, 1976.

\bibitem{JEB1977Prikry}
James~E. Baumgartner and Karel Prikry.
\newblock Singular cardinals and the generalized continuum hypothesis.
\newblock {\em Amer. Math. Monthly}, 84(2):108--113, 1977.

\bibitem{JEB1987Shelah}
James~E. Baumgartner and Saharon Shelah.
\newblock Remarks on superatomic {B}oolean algebras.
\newblock {\em Ann. Pure Appl. Logic}, 33(2):109--129, 1987.

\bibitem{JEB1993ST}
James~E. Baumgartner, Saharon Shelah, and Simon Thomas.
\newblock Maximal subgroups of infinite symmetric groups.
\newblock {\em Notre Dame J. Formal Logic}, 34(1):1--11, 1993.

\bibitem{JEB1991Spinas}
James~E. Baumgartner and Otmar Spinas.
\newblock Independence and consistency proofs in quadratic form theory.
\newblock {\em J. Symbolic Logic}, 56(4):1195--1211, 1991.

\bibitem{JEB2002Tall}
James~E. Baumgartner and Franklin~D. Tall.
\newblock Reflecting {L}indel\"ofness.
\newblock In {\em Proceedings of the {I}nternational {C}onference on {T}opology
  and its {A}pplications ({Y}okohama, 1999)}, volume 122, pages 35--49, 2002.

\bibitem{JEB1978Taylor}
James~E. Baumgartner and Alan~D. Taylor.
\newblock Partition theorems and ultrafilters.
\newblock {\em Trans. Amer. Math. Soc.}, 241:283--309, 1978.

\bibitem{JEB1982Tsatideals1}
James~E. Baumgartner and Alan~D. Taylor.
\newblock Saturation properties of ideals in generic extensions. {I}.
\newblock {\em Trans. Amer. Math. Soc.}, 270(2):557--574, 1982.

\bibitem{JEB1982Tsatideals2}
James~E. Baumgartner and Alan~D. Taylor.
\newblock Saturation properties of ideals in generic extensions. {II}.
\newblock {\em Trans. Amer. Math. Soc.}, 271(2):587--609, 1982.

\bibitem{JEB1977TW}
James~E. Baumgartner, Alan~D. Taylor, and Stanley Wagon.
\newblock On splitting stationary subsets of large cardinals.
\newblock {\em J. Symbolic Logic}, 42(2):203--214, 1977.

\bibitem{JEB1978TW}
James~E. Baumgartner, Alan~D. Taylor, and Stanley Wagon.
\newblock Ideals on uncountable cardinals.
\newblock In {\em Logic {C}olloquium '77 ({P}roc. {C}onf., {W}roc\l aw, 1977)},
  volume~96 of {\em Stud. Logic Foundations Math.}, pages 67--77.
  North-Holland, Amsterdam, 1978.

\bibitem{JEB1990vD}
James~E. Baumgartner and Eric~K. van Douwen.
\newblock Strong realcompactness and weakly measurable cardinals.
\newblock {\em Topology Appl.}, 35(2-3):239--251, 1990.

\bibitem{JEB1982Weese}
James~E. Baumgartner and Martin Weese.
\newblock Partition algebras for almost-disjoint families.
\newblock {\em Trans. Amer. Math. Soc.}, 274(2):619--630, 1982.

\bibitem{BellerLitman1980}
Aaron Beller and Ami Litman.
\newblock A strengthening of {J}ensen's {$ cm$} principles.
\newblock {\em J. Symbolic Logic}, 45(2):251--264, 1980.

\bibitem{BrechKoszmider}
Christina Brech and Piotr Koszmider.
\newblock On universal {B}anach spaces of density continuum.
\newblock {\em Israel J. Math.}, 190:93--110, 2012.

\bibitem{Brodsky2014}
Ari~M. Brodsky.
\newblock A theory of stationary trees and the balanced
  {B}aumgartner-{H}ajnal-{T}odorcevic theorem for trees.
\newblock {\em Acta Math. Hungar.}, 144(2):285--352, 2014.

\bibitem{ChangKeisler1973}
Chen~Chung Chang and H.~Jerome Keisler.
\newblock {\em Model theory}.
\newblock North-Holland Publishing Co., Amsterdam-London; American Elsevier
  Publishing Co., Inc., New York, 1973.
\newblock Studies in Logic and the Foundations of Mathematics, Vol. 73.

\bibitem{ClaverieSchindler2012}
Benjamin Claverie and Ralf Schindler.
\newblock Woodin's axiom {$(\ast)$}, bounded forcing axioms, and precipitous
  ideals on {$\omega_1$}.
\newblock {\em J. Symbolic Logic}, 77(2):475--498, 2012.

\bibitem{Cummings2010hand}
James Cummings.
\newblock Iterated forcing and elementary embeddings.
\newblock In {\em Handbook of set theory}, pages 775--883. Springer, Dordrecht,
  2010.

\bibitem{Devlin1983Yguide}
Keith~J. Devlin.
\newblock The {Y}orkshireman's guide to proper forcing.
\newblock In {\em Surveys in set theory}, volume~87 of {\em London Math. Soc.
  Lecture Note Ser.}, pages 60--115. Cambridge Univ. Press, Cambridge, 1983.

\bibitem{DevlinJohns1974}
Keith~J. Devlin and H{\aa}vard Johnsbr{\aa}ten.
\newblock {\em The {S}ouslin problem}.
\newblock Lecture Notes in Mathematics, Vol. 405. Springer-Verlag, Berlin-New
  York, 1974.

\bibitem{Easton1964}
William~B. Easton.
\newblock {\em Powers of Regular Cardinals}.
\newblock PhD thesis, Princeton University, Princeton, N.J., 1964.
\newblock Advisor {A}.~Church.

\bibitem{Erdos42}
Paul Erd\H{o}s.
\newblock Some set-theoretical properties of graphs.
\newblock {\em Revista de la Universidad Nacional de Tucum{{\'a}}n, Serie A},
  3:363--367, 1942.

\bibitem{EH74.prob}
Paul Erd\H{o}s and And{\'{a}}s Hajnal.
\newblock Solved and unsolved problems in set theory.
\newblock In {\em Proceedings of the Tarski Symposium, (Proc. Sympos. Pure
  Math., UC Berkeley, CA, 1971)}, volume~25, pages 261--265. Amer. Math. Soc.,
  Providence, R.I., 1974.

\bibitem{EH71.prob}
Paul Erd\H{o}s and Andr{\'{a}}s Hajnal.
\newblock Unsolved problems in set theory.
\newblock In {\em Axiomatic {S}et {T}heory ({P}roc. {S}ympos. {P}ure {M}ath.,
  {UCLA}, 1967)}, volume~13, pages 17--48. Amer. Math. Soc., Providence, R.I.,
  1971.

\bibitem{EP1965HajnalRado}
Paul Erd\H{o}s, Andr{\'{a}}s Hajnal, and Richard Rado.
\newblock Partition relations for cardinal numbers.
\newblock {\em Acta Math. Acad. Sci. Hungar.}, 16:93--196, 1965.

\bibitem{ErdosRado1953}
Paul Erd\H{o}s and Richard Rado.
\newblock A problem on ordered sets.
\newblock {\em J. London Math. Soc. (2)}, 28(4):426--438, 1953.

\bibitem{EP1956Rado}
Paul Erd\H{o}s and Richard Rado.
\newblock A partition calculus in set theory.
\newblock {\em Bull. Amer. Math. Soc.}, 62:427--489, 1956.

\bibitem{ErdosTarski1943}
Paul Erd\H{o}s and Alfred Tarski.
\newblock On families of mutually exclusive sets.
\newblock {\em Ann. of Math. (2)}, 44(2):315--329, April 1943.

\bibitem{EP1975combo.prob}
Paul Erd{\H{o}}s.
\newblock Problems and results on finite and infinite combinatorial analysis.
\newblock In {\em Infinite and finite sets ({C}olloq., {K}eszthely, 1973;
  dedicated to {P}. {E}rd{\H o}s on his 60th birthday), {V}ol. {I}}, pages
  403--424. Colloq. Math. Soc. J\'anos Bolyai, Vol. 10. North-Holland,
  Amsterdam, 1975.

\bibitem{EP1961Hajnal}
Paul Erd{\H{o}}s and Andr{\'{a}}s Hajnal.
\newblock On a property of families of sets.
\newblock {\em Acta Math. Acad. Sci. Hungar}, 12:87--123, 1961.

\bibitem{ErdosHajnal1962Ramsey}
Paul Erd{\H{o}}s and Andr{\'a}s Hajnal.
\newblock Some remarks concerning our paper ``{O}n the structure of
  set-mappings''. {N}on-existence of a two-valued {$\sigma $}-measure for the
  first uncountable inaccessible cardinal.
\newblock {\em Acta Math. Acad. Sci. Hungar.}, 13:223--226, 1962.

\bibitem{EP1966Hajnal}
Paul Erd{\H{o}}s and Andr{\'{a}}s Hajnal.
\newblock On chromatic number of graphs and set-systems.
\newblock {\em Acta Math. Acad. Sci. Hungar}, 17:61--99, 1966.

\bibitem{ErdosRado1956}
Paul Erd{\"o}s and Richard Rado.
\newblock A partition calculus in set theory.
\newblock {\em Bull. Amer. Math. Soc.}, 62:427--489, 1956.

\bibitem{ErdosTarski1961}
Paul Erd{\H{o}}s and Alfred Tarski.
\newblock On some problems involving inaccessible cardinals.
\newblock In {\em Essays on the foundations of mathematics}, pages 50--82.
  Magnes Press, Hebrew Univ., Jerusalem, 1961.

\bibitem{Eskew2014PhDthesis}
Monroe Eskew.
\newblock {\em Measurability {P}roperties on {S}mall {C}ardinals}.
\newblock PhD thesis, University of California, Irvine, CA, 2014.
\newblock Advisor {M}.~{Z}eman.

\bibitem{Feng1990}
Qi~Feng.
\newblock A hierarchy of {R}amsey cardinals.
\newblock {\em Ann. Pure Appl. Logic}, 49(3):257--277, 1990.

\bibitem{FleissnerKunen1978}
William~G. Fleissner and Kenneth Kunen.
\newblock Barely {B}aire spaces.
\newblock {\em Fund. Math.}, 101(3):229--240, 1978.

\bibitem{Foreman2015}
Matthew Foreman.
\newblock Chang's conjecture, generic elementary embeddings and inner models
  for huge cardinals.
\newblock {\em Bull. Symb. Log.}, 21(3):251--269, 2015.

\bibitem{Foreman2003Hajnal}
Matthew Foreman and Andr{\'{a}}s Hajnal.
\newblock A partition relation for successors of large cardinals.
\newblock {\em Math. Ann.}, 325(3):583--623, 2003.

\bibitem{ForemanLaver1988}
Matthew Foreman and Richard Laver.
\newblock Some downwards transfer properties for {$\aleph_2$}.
\newblock {\em Adv. in Math.}, 67(2):230--238, 1988.

\bibitem{Foreman1995Magidor}
Matthew Foreman and Menachem Magidor.
\newblock Large cardinals and definable counterexamples to the continuum
  hypothesis.
\newblock {\em Ann. Pure Appl. Logic}, 76(1):47--97, 1995.

\bibitem{Foreman1988MagShelah}
Matthew Foreman, Menachem Magidor, and Saharon Shelah.
\newblock Martin's maximum, saturated ideals and nonregular ultrafilters. {II}.
\newblock {\em Ann. of Math. (2)}, 127(3):521--545, 1988.

\bibitem{FriedmanH1974}
Harvey~M. Friedman.
\newblock On closed sets of ordinals.
\newblock {\em Proc. Amer. Math. Soc.}, 43:190--192, 1974.

\bibitem{FriedmanH2001}
Harvey~M. Friedman.
\newblock Subtle cardinals and linear orderings.
\newblock {\em Ann. Pure Appl. Logic}, 107(1-3):1--34, 2001.

\bibitem{Friedman2006}
Sy-David Friedman.
\newblock Forcing with finite conditions.
\newblock In {\em Set theory}, Trends Math., pages 285--295. Birkh\"auser,
  Basel, 2006.

\bibitem{GaifmanSpecker1964}
Haim Gaifman and Ernst Specker.
\newblock Isomorphism types of trees.
\newblock {\em Proc. Amer. Math. Soc.}, 15:1--7, 1964.

\bibitem{Galvin1973chrom}
Fred Galvin.
\newblock Chromatic numbers of subgraphs.
\newblock {\em Period. Math. Hungar.}, 4:117--119, 1973.

\bibitem{Galvin1975BH}
Fred Galvin.
\newblock On a partition theorem of {B}aumgartner and {H}ajnal.
\newblock In {\em Infinite and {F}inite {S}ets ({C}olloq., {K}eszthely, 1973;
  dedicated to {P}. {E}rd{\H o}s on his 60th birthday), {V}ol. {II}}, pages
  711--729. Colloq. Math. Soc. J{\'a}nos Bolyai, Vol. 10. North-Holland,
  Amsterdam, 1975.

\bibitem{Gitik2010normal}
Moti Gitik.
\newblock On normal precipitous ideals.
\newblock {\em Israel J. Math.}, 175:191--219, 2010.

\bibitem{Gitik1997Shelah}
Moti Gitik and Saharon Shelah.
\newblock Less saturated ideals.
\newblock {\em Proc. Amer. Math. Soc.}, 125(5):1523--1530, 1997.

\bibitem{Hajnal2013}
Andr\'as Hajnal.
\newblock Private conversation on {B}aumgartner-{H}ajnal {T}heorem.
\newblock Jean Larson interviewed Hajnal at the Erd\H{o}s Centennial in
  Budapest., July 2013.

\bibitem{Halbeisen2005}
Lorenz Halbeisen.
\newblock Families of almost disjoint {H}amel bases.
\newblock {\em Extracta Math.}, 20(2):199--202, 2005.

\bibitem{HarringtonShelah1985}
Leo Harrington and Saharon Shelah.
\newblock Some exact equiconsistency results in set theory.
\newblock {\em Notre Dame J. Formal Logic}, 26(2):178--188, 1985.

\bibitem{Hausdorff1907a}
Felix Hausdorff.
\newblock {U}ntersuchungen {{\"u}}ber {O}rdnugstypen {IV}, {V}.
\newblock {\em {B}erichte {\"u}ber die {V}erlhandlungen der {K}{{\"o}}niglich
  {S}{{\"a}}sischen {G}esellschaft der {W}issenschaften zu {L}eipzig,
  {M}athematisch-{P}hysiche {K}lasse}, 59:84--159, 1907.

\bibitem{Hausdorff1908}
Felix Hausdorff.
\newblock Grundz{{\"u}}ge einer {T}heorie der geordnete {M}engenlehre.
\newblock {\em Math. Ann.}, 65:435--505, 1908.

\bibitem{Ishiu2005AxiomA}
Tetsuya Ishiu.
\newblock {$\alpha$}-properness and {A}xiom {A}.
\newblock {\em Fund. Math.}, 186(1):25--37, 2005.

\bibitem{Jech1967Suslin}
Thomas Jech.
\newblock Non-provability of {S}ouslin's hypothesis.
\newblock {\em Comment. Math. Univ. Carolinae}, 8(2):291--305, 1967.

\bibitem{Jech2003}
Thomas Jech.
\newblock {\em Set theory}.
\newblock Springer Monographs in Mathematics. Springer-Verlag, Berlin, 2003.
\newblock The third millennium edition, revised and expanded.

\bibitem{JechPrikry1976}
Thomas Jech and Karel Prikry.
\newblock On ideals of sets and the power set operation.
\newblock {\em Bull. Amer. Math. Soc.}, 82(4):593--595, 1976.

\bibitem{UCLAvol2}
Thomas~J. Jech, editor.
\newblock {\em Axiomatic set theory}.
\newblock {P}roc. {S}ympos. {P}ure {M}ath., {V}ol. {XIII}, {P}art {II}, {U}niv.
  {C}alifornia, {L}os {A}ngeles, {C}alif., 1967. American Mathematical Society,
  Providence, R.I., 1974.

\bibitem{Jensen1967book}
Ronald~B. Jensen.
\newblock {\em Modelle der {M}engenlehre. {W}iderspruchsfreiheit und
  {U}nabh\"angigkeit der {K}ontinuum-{H}ypothese und des {A}uswahlaxioms}.
\newblock Ausgearbeitet von Franz Josef Leven. Lecture Notes in Mathematics,
  No. 37. Springer-Verlag, Berlin, 1967.

\bibitem{Jensen1972fine}
Ronald~B. Jensen.
\newblock The fine structure of the constructible hierarchy.
\newblock {\em Ann. Math. Logic}, 4:229--308; erratum, ibid. 4 (1972), 443,
  1972.
\newblock With a section by Jack Silver.

\bibitem{JensenKunen1969}
Ronald~B. Jensen and Kenneth Kunen.
\newblock Some {C}ombinatorial {P}roperties of {$L$} and {$V$}.
\newblock Handwritten manuscript, scanned and posted to the Jensen website,
  April 1969.

\bibitem{StackingMice2009}
Ronald~B. Jensen, Ernest Schimmerling, Ralf Schindler, and John~R. Steel.
\newblock Stacking mice.
\newblock {\em J. Symbolic Logic}, 74(1):315--335, 2009.

\bibitem{JensenSteel2013}
Ronald~B. Jensen and John~R. Steel.
\newblock {$K$} without the measurable.
\newblock {\em J. Symbolic Logic}, 78(3):708--734, 2013.

\bibitem{Jones2006}
Albin~L. Jones.
\newblock A polarized partition relation for weakly compact cardinals using
  elementary substructures.
\newblock {\em J. Symbolic Logic}, 71(4):1342--1352, 2006.

\bibitem{Kanamori2003}
Akihiro Kanamori.
\newblock {\em The higher infinite}.
\newblock Springer Monographs in Mathematics. Springer-Verlag, Berlin, second
  edition, 2009.
\newblock Large cardinals in set theory from their beginnings, Paperback
  reprint of the 2003 edition.

\bibitem{Kanamori1978Magidor}
Akihiro Kanamori and Menachem Magidor.
\newblock The evolution of large cardinal axioms in set theory.
\newblock In {\em Higher set theory ({P}roc. {C}onf., {M}ath.
  {F}orschungsinst., {O}berwolfach, 1977)}, volume 669 of {\em Lecture Notes in
  Math.}, pages 99--275. Springer, Berlin, 1978.

\bibitem{KeislerTarski1963}
H.~J. Keisler and A.~Tarski.
\newblock From accessible to inaccessible cardinals. {R}esults holding for all
  accessible cardinal numbers and the problem of their extension to
  inaccessible ones.
\newblock {\em Fund. Math.}, 53:225--308, 1963/1964.
\newblock Corrections on page 119 of volume 55, 1965.

\bibitem{Komjath1988}
P{\'e}ter Komj{\'a}th.
\newblock Consistency results on infinite graphs.
\newblock {\em Israel J. Math.}, 61, 1988.

\bibitem{Komjath2002}
P{\'e}ter Komj{\'a}th.
\newblock Subgraph chromatic number.
\newblock In {\em Set theory ({P}iscataway, {NJ}, 1999)}, volume~58 of {\em
  DIMACS Ser. Discrete Math. Theoret. Comput. Sci.}, pages 99--106. Amer. Math.
  Soc., Providence, RI, 2002.

\bibitem{Krueger2003}
John Krueger.
\newblock Fat sets and saturated ideals.
\newblock {\em J. Symbolic Logic}, 68(3):837--845, 2003.

\bibitem{Krueger2015clubII}
John Krueger.
\newblock Adding a club with finite conditions, {P}art {II}.
\newblock {\em Arch. Math. Logic}, 54(1-2):161--172, 2015.

\bibitem{KruegerSchimmerling2011}
John Krueger and Ernest Schimmerling.
\newblock An equiconsistency result on partial squares.
\newblock {\em J. Math. Log.}, 11(1):29--59, 2011.

\bibitem{Kunen1970}
Kenneth Kunen.
\newblock Some applications of iterated ultrapowers in set theory.
\newblock {\em Ann. Math. Logic}, 1:179--227, 1970.

\bibitem{Kunen1978}
Kenneth Kunen.
\newblock Saturated ideals.
\newblock {\em J. Symbolic Logic}, 43(1):65--76, 1978.

\bibitem{Kurepa1935}
{\Dbar}uro Kurepa.
\newblock Ensembles ordonn{\'e}es et ramifi{\'e}s.
\newblock {\em Publ. Math. Univ. Belgrade}, 4:1--138, 1935.
\newblock A35.

\bibitem{Kurepa1937Atree}
{\Dbar}uro Kurepa.
\newblock Ensembles lin{{\'e}}aires et une classe de tableaux ramifi{{\'e}}s
  ({T}ableaux ramifi{{\'e}}s de {M}. {A}ronszajn).
\newblock {\em Publ. Inst. Math. (Beograd)}, 6:129--160, 1937.

\bibitem{Kurepa1939V}
{\Dbar}uro Kurepa.
\newblock Sur la puissance des ensembles partillement ordonn{\'e}s.
\newblock {\em C. R. Soc. Sci. Varsovie, Cl. Math.}, 32:61--67, 1939.
\newblock Sometimes the journal is listed in {P}olish: {S}prawozdania
  {T}owarzystwo {N}aukowe {W}arszawa {M}at.-{F}iz. as in the Math Review by
  Bagemihl.

\bibitem{Larson2000}
Paul Larson.
\newblock Separating stationary reflection principles.
\newblock {\em J. Symbolic Logic}, 65(1):247--258, 2000.

\bibitem{Laverthesis1969}
Richard Laver.
\newblock {\em Order {T}ypes and {W}ell-{Q}uasi-{O}rderings}.
\newblock PhD thesis, University of California, Berkeley, 1969.
\newblock Advisor {R}.~{M}c{K}enzie.

\bibitem{Laver1976BorelConj}
Richard Laver.
\newblock On the consistency of {B}orel's conjecture.
\newblock {\em Acta Math.}, 137(3-4):151--169, 1976.

\bibitem{Laver1978diamond}
Richard Laver.
\newblock Making the supercompactness of {$\kappa $} indestructible under
  {$\kappa $}-directed closed forcing.
\newblock {\em Israel J. Math.}, 29(4):385--388, 1978.

\bibitem{Magidor1971a}
Menachem Magidor.
\newblock There are many normal ultrafiltres corresponding to a supercompact
  cardinal.
\newblock {\em Israel J. Math.}, 9:186--192, 1971.

\bibitem{Magidor1982}
Menachem Magidor.
\newblock Reflecting stationary sets.
\newblock {\em J. Symbolic Logic}, 47(4):755--771 (1983), 1982.

\bibitem{Malitz1968}
Jerome~I. Malitz.
\newblock The {H}anf number for complete ${L}_{\omega_1,\omega}$-sentences.
\newblock In {\em The {S}yntax and {S}emantics of {I}nfinitary {L}anguages},
  volume~72 of {\em Lecture Notes in Mathematics}, pages 166--181.
  Springer-Verlag, Berlin, 1968.

\bibitem{Martinez2001}
Juan~Carlos Mart{\'{\i}}nez.
\newblock A consistency result on thin-very tall {B}oolean algebras.
\newblock {\em Israel J. Math.}, 123:273--284, 2001.

\bibitem{Martinez2016}
Juan~Carlos Mart{\'{\i}}nez.
\newblock On cardinal sequences of {LCS} spaces.
\newblock {\em Topology Appl.}, 203:91--97, 2016.

\bibitem{Matet2003}
Pierre Matet.
\newblock Partition relations for {$\kappa$}-normal ideals on
  {$P_\kappa(\lambda)$}.
\newblock {\em Ann. Pure Appl. Logic}, 121(1):89--111, 2003.

\bibitem{MilnerPrikry1991}
Eric~C. Milner and Karel Prikry.
\newblock A partition relation for triples using a model of {T}odor\v cevi\'c.
\newblock {\em Discrete Math.}, 95(1-3):183--191, 1991.
\newblock Directions in infinite graph theory and combinatorics (Cambridge,
  1989).

\bibitem{Mitchell1970}
William~J. Mitchell.
\newblock {\em Aronszajn trees and the independence of the transfer property}.
\newblock PhD thesis, University of California, Berkeley, 1970.
\newblock Advisor {J}.~{S}ilver.

\bibitem{MitchellAtree1972}
William~J. Mitchell.
\newblock Aronszajn trees and the independence of the transfer property.
\newblock {\em Ann. Math. Logic}, 5:21--46, 1972.

\bibitem{Mitchell1972Atree}
William~J. Mitchell.
\newblock Aronszajn trees and the independence of the transfer property.
\newblock {\em Ann. Math. Logic}, 5:21--46, 1972/73.

\bibitem{Mitchell2009}
William~J. Mitchell.
\newblock {$I[\omega_2]$} can be the nonstationary ideal on {${\rm
  Cof}(\omega_1)$}.
\newblock {\em Trans. Amer. Math. Soc.}, 361(2):561--601, 2009.

\bibitem{Monk2006}
J.~Donald Monk.
\newblock The size of maximal almost disjoint families.
\newblock {\em Dissertationes Math. (Rozprawy Mat.)}, 437:47, 2006.

\bibitem{MooreG1988}
Gregory~H. Moore.
\newblock The origins of forcing.
\newblock In {\em Logic {C}olloquium '86 ({H}ull 1986)}, volume 124 of {\em
  Studies in Logic and the Foundations of Mathematics}, pages 143--173. North
  Holland, Amsterdam, 1988.
\newblock {F}rank~{R}. {D}rake and {J}ohn~{K}. {T}russ, editors.

\bibitem{MooreJT2006Basis}
Justin~Tatch Moore.
\newblock A five element basis for the uncountable linear orders.
\newblock {\em Ann. of Math. (2)}, 163(2):669--688, 2006.

\bibitem{Neeman2014}
Itay Neeman.
\newblock Forcing with sequences of models of two types.
\newblock {\em Notre Dame J. Form. Log.}, 55(2):265--298, 2014.

\bibitem{Ramsey1930}
Frank~P. Ramsey.
\newblock On a {P}roblem of {F}ormal {L}ogic.
\newblock {\em Proc. London Math. Soc.}, S2-30(1):264--286, 1930.

\bibitem{Rinot2014}
Assaf Rinot.
\newblock Chain conditions of products, and weakly compact cardinals.
\newblock {\em Bull. Symb. Log.}, 20(3):293--314, 2014.

\bibitem{Rinot2015}
Assaf Rinot.
\newblock Chromatic numbers of graphs---large gaps.
\newblock {\em Combinatorica}, 35(2):215--233, 2015.

\bibitem{UCLAvol1}
Dana~J. Scott, editor.
\newblock {\em Axiomatic set theory}.
\newblock {P}roc. {S}ympos. {P}ure {M}ath {V}ol. {XIII}, {P}art {I}, {U}niv.
  {C}alif., {L}os {A}ngeles, Calif., 1967. American Mathematical Society,
  Providence, R.I., 1971.

\bibitem{SharpeWelch2011}
Ian Sharpe and Philip~D. Welch.
\newblock Greatly {E}rd{\H o}s cardinals with some generalizations to the
  {C}hang and {R}amsey properties.
\newblock {\em Ann. Pure Appl. Logic}, 162(11):863--902, 2011.

\bibitem{Shelah1980LinO}
Saharon Shelah.
\newblock Independence results.
\newblock {\em J. Symbolic Logic}, 45(3):563--573, 1980.

\bibitem{Shelah1982ProperForcing}
Saharon Shelah.
\newblock {\em Proper forcing}, volume 940 of {\em Lecture Notes in
  Mathematics}.
\newblock Springer-Verlag, Berlin, 1982.

\bibitem{Shelah1987NSsat}
Saharon Shelah.
\newblock Iterated forcing and normal ideals on {$\omega_1$}.
\newblock {\em Israel J. Math.}, 60(3):345--380, 1987.

\bibitem{Shelah1991partial}
Saharon Shelah.
\newblock Reflecting stationary sets and successors of singular cardinals.
\newblock {\em Arch. Math. Logic}, 31(1):25--53, 1991.

\bibitem{ShelahProperImproper}
Saharon Shelah.
\newblock {\em Proper and improper forcing}.
\newblock Perspectives in Mathematical Logic. Springer-Verlag, Berlin, second
  edition, 1998.

\bibitem{ShelahWoodin1990}
Saharon Shelah and Hugh Woodin.
\newblock Large cardinals imply that every reasonably definable set of reals is
  {L}ebesgue measurable.
\newblock {\em Israel J. Math.}, 70(3):381--394, 1990.

\bibitem{Shoenfield1971}
Joseph~R. Shoenfield.
\newblock Unramified forcing.
\newblock In {\em Axiomatic {S}et {T}heory ({P}roc. {S}ympos. {P}ure {M}ath.,
  {V}ol. {XIII}, {P}art {I}, {U}niv. {C}alifornia, {L}os {A}ngeles, {C}alif.,
  1967)}, pages 357--381. Amer. Math. Soc., Providence, R.I., 1971.

\bibitem{Sierpinski1928ad}
Wac{\l}aw Sierpi{{\'n}}ski.
\newblock Sur une d{\'e}composition d'ensembles.
\newblock {\em Monatsh. Math. Phys.}, 35(1):239--242, 1928.

\bibitem{Sierpinski1933}
Wac{\l}aw Sierpi{\'n}ski.
\newblock Sur un probl\`eme de la th\'eorie des relations.
\newblock {\em Ann. Scuola Norm. Sup. Pisa Cl. Sci. (2)}, 2(3):285--287, 1933.

\bibitem{Sierpinski1956}
Wac{\l}aw Sierpi{\'n}ski.
\newblock {\em Hypoth\`ese du continu}.
\newblock Chelsea Publishing Company, New York, N. Y., 1956.
\newblock 2nd ed; 1st ed,1934, Warsaw.

\bibitem{Silver1966PhDthesis}
Jack Silver.
\newblock {\em Some applications of model theory in set theory}.
\newblock PhD thesis, University of California, Berkeley, 1966.
\newblock Advisor {R}.~{V}aught.

\bibitem{Silver1970}
Jack Silver.
\newblock A large cardinal in the constructible universe.
\newblock {\em Fund. Math.}, 69:93--100, 1970.

\bibitem{Silver1975}
Jack Silver.
\newblock On the singular cardinals problem.
\newblock In {\em Proceedings of the {I}nternational {C}ongress of
  {M}athematicians ({V}ancouver, {B}. {C}., 1974), {V}ol. 1}, pages 265--268.
  Canad. Math. Congress, Montreal, Que., 1975.

\bibitem{Solovay1964Julabs}
Robert~M. Solovay.
\newblock The measure problem.
\newblock {\em J. Symbolic Logic}, 29:227--228, 1964.
\newblock Abstract for the {A}ssociation of {S}ymbolic {L}ogic meeting of July
  13--17, 1964 held at the University of Bristol; received July 6, 1964.

\bibitem{Solovay1970}
Robert~M. Solovay.
\newblock A model of set-theory in which every set of reals is {L}ebesgue
  measurable.
\newblock {\em Ann. of Math. (2)}, 92:1--56, 1970.

\bibitem{Solovay1971}
Robert~M. Solovay.
\newblock Real-valued measurable cardinals.
\newblock In {\em Axiomatic {S}et {T}heory ({P}roc. {S}ympos. {P}ure {M}ath.,
  {V}ol. XIII, Part I, Univ. California, Los Angeles, Calif., 1967)}, pages
  397--428. Amer. Math. Soc., Providence, R.I., 1971.

\bibitem{ReinhardtSolovayKanamori1978}
Robert~M. Solovay, William~N. Reinhardt, and Akihiro Kanamori.
\newblock Strong axioms of infinity and elementary embeddings.
\newblock {\em Ann. Math. Logic}, 13(1):73--116, 1978.

\bibitem{SolovayTennenbaum1971}
Robert~M. Solovay and Stanley Tennenbaum.
\newblock Iterated {C}ohen extensions and {S}ouslin's problem.
\newblock {\em Ann. of Math. (2)}, 94:201--245, 1971.

\bibitem{Steel1996CoreIterate}
John~R. Steel.
\newblock {\em The core model iterability problem}, volume~8 of {\em Lecture
  Notes in Logic}.
\newblock Springer-Verlag, Berlin, 1996.

\bibitem{Steel2007}
John~R. Steel.
\newblock What is {$\dots$} a {W}oodin cardinal?
\newblock {\em Notices Amer. Math. Soc.}, 54(9):1146--1147, 2007.

\bibitem{Suslin1920}
Mikhail Suslin.
\newblock Probl{\`e}me 3.
\newblock {\em Fund. Math.}, 1, 1920.

\bibitem{Tarski1925gch}
Alfred Tarski.
\newblock Quelques th{{\'e}}or{\`e}mes qui {\'e}quivalent {\`a} l'axiom du
  choix.
\newblock {\em Fund. Math.}, 7:147--154, 1925.

\bibitem{Tarski1928ad}
Alfred Tarski.
\newblock Sur la d{{\'e}}composition des ensembles en sous-ensembles presque
  disjoints.
\newblock {\em Fund. Math.}, 12:188--205, 1928.

\bibitem{Tarski1929ad}
Alfred Tarski.
\newblock Sur la d{{\'e}}composition des ensembles en sous-ensembles presque
  disjoints (suppl\'ement).
\newblock {\em Fund. Math.}, 14:205--215, 1929.

\bibitem{Tarski1945}
Alfred Tarski.
\newblock Ideale in vollst{\"a}digen {M}engenk{\"o}rpern, ii.
\newblock {\em Fund. Math.}, 33:51--65, 1945.

\bibitem{Tennenbaum1968Suslin}
Stanley Tennenbaum.
\newblock Souslin's problem.
\newblock {\em Proc. Nat. Acad. Sci. U.S.A.}, 59:60--63, 1968.

\bibitem{Todor1981wKH}
Stevo Todorcevic.
\newblock Some consequences of {${\rm MA}+\neg{\rm wKH}$}.
\newblock {\em Topology Appl.}, 12(2):187--202, 1981.

\bibitem{Todor1983forcePos}
Stevo Todorcevic.
\newblock Forcing positive partition relations.
\newblock {\em Trans. Amer. Math. Soc.}, 280(2):703--720, 1983.

\bibitem{Todor1984onPFA}
Stevo Todorcevic.
\newblock A note on the proper forcing axiom.
\newblock In James~E. Baumgartner, Donald~A. Martin, and Saharon Shelah,
  editors, {\em Axiomatic set theory}, volume~31 of {\em Contemporary
  Mathematics}, pages 209--218. American Mathematical Society, 1984.

\bibitem{Todor1984settop}
Stevo Todorcevic.
\newblock Trees and linearly ordered sets.
\newblock In {\em Handbook of set-theoretic topology}, pages 235--293.
  North-Holland, Amsterdam, 1984.

\bibitem{Todor1985poset}
Stevo Todorcevic.
\newblock Partition relations for partially ordered sets.
\newblock {\em Acta Math.}, 155(1-2):1--25, 1985.

\bibitem{Todor1987partPairs}
Stevo Todorcevic.
\newblock Partitioning pairs of countable ordinals.
\newblock {\em Acta Math.}, 159(3-4):261--294, 1987.

\bibitem{TodorProbTop1989}
Stevo Todorcevic.
\newblock {\em Partition {P}roblems in {T}opology}, volume~84 of {\em
  Contemporary Mathematics}.
\newblock American Mathematical Society, Providence, RI, 1989.

\bibitem{Todor1997}
Stevo Todorcevic.
\newblock Comparing the continuum with the first two uncountable cardinals.
\newblock In {\em Logic and scientific methods ({F}lorence, 1995)}, volume 259
  of {\em Synthese Lib.}, pages 145--155. Kluwer Acad. Publ., Dordrecht, 1997.

\bibitem{Ulam1930}
Stanislaw Ulam.
\newblock Zur masstheorie in der allgemeinen mengenlehre.
\newblock {\em Fund. Math.}, 16(1):140--150, 1930.

\bibitem{Velickovic1992forcestat}
Boban Veli{\v{c}}kovi{\'c}.
\newblock Forcing axioms and stationary sets.
\newblock {\em Adv. Math.}, 94(2):256--284, 1992.

\bibitem{VelickovicVenturi}
Boban Velickovic and Giorgio Venturi.
\newblock Proper forcing remastered.
\newblock In {\em Appalachian {S}et {T}heory 2006-2012}, volume 213 of {\em
  London Mathematical Society Lecture Note Series}, pages 331--362. Cambridge
  University Press, Cambridge, 2013.

\bibitem{Viale2011Weiss}
Matteo Viale and Christoph Wei{\ss}.
\newblock On the consistency strength of the proper forcing axiom.
\newblock {\em Adv. Math.}, 228(5):2672--2687, 2011.

\end{thebibliography}
\end{document}